\documentclass[11pt,a4paper]{article}
\usepackage[a4paper]{geometry}
\usepackage{amssymb,latexsym,amsmath,amsfonts,amsthm}
\usepackage{graphicx}
\usepackage{dsfont}
\usepackage{mathrsfs}
\usepackage{epsfig}
\usepackage{epstopdf}
\usepackage{booktabs}
\usepackage{mathtools}
\usepackage{comment,verbatim}
\usepackage{overpic}
\usepackage{hyperref}
\usepackage{bbm}
\usepackage[toc,page]{appendix}

\newtheorem{theorem}{Theorem}[section]
\newtheorem{lemma}[theorem]{Lemma}

\theoremstyle{definition}

\theoremstyle{definition}

\theoremstyle{remark}

\newtheorem{remark}[theorem]{Remark}

\numberwithin{equation}{section}

\usepackage{color}

\usepackage[normalem]{ulem}

\def\Xint#1{\mathchoice
{\XXint\displaystyle\textstyle{#1}}%
{\XXint\textstyle\scriptstyle{#1}}%
{\XXint\scriptstyle
\scriptscriptstyle{#1}}%
{\XXint\scriptscriptstyle
\scriptscriptstyle{#1}}%
\!\int}
\def\XXint#1#2#3{{
\setbox0=\hbox{$#1{#2#3}{\int}$}
\vcenter{\hbox{$#2#3$}}\kern-.5\wd0}}
\def\dashint{\Xint-}

\begin{document}
\title{Convergence analysis of Laguerre approximations for analytic functions}
\author{Haiyong Wang\footnotemark[1]~\footnotemark[2]}
\date{}

\maketitle

\footnotetext[1]{School of Mathematics and Statistics, Huazhong
University of Science and Technology, Wuhan 430074, P. R. China.
E-mail: \texttt{haiyongwang@hust.edu.cn}}

\footnotetext[2]{Hubei Key Laboratory of Engineering Modeling and
Scientific Computing, Huazhong University of Science and Technology,
Wuhan 430074, China}

\begin{abstract}
Laguerre spectral approximations play an important role in the
development of efficient algorithms for problems in unbounded
domains. In this paper, we present a comprehensive convergence rate
analysis of Laguerre spectral approximations for analytic functions.
By exploiting contour integral techniques from complex analysis, we
prove that Laguerre projection and interpolation methods of degree
$n$ converge at the root-exponential rate $O(\exp(-2\rho\sqrt{n}))$
with $\rho>0$ when the underlying function is analytic inside and on
a parabola with focus at the origin and vertex at $z=-\rho^2$. As
far as we know, this is the first rigorous proof of root-exponential
convergence of Laguerre approximations for analytic functions.
Several important applications of our analysis are also discussed,
including Laguerre spectral differentiations, Gauss-Laguerre
quadrature rules, the scaling factor and the Weeks method for the inversion of Laplace
transform, and some sharp convergence rate estimates are derived.
Numerical experiments are presented to verify the theoretical
results.
\end{abstract}


\vspace{0.05in}


\section{Introduction}\label{sec:introduction}
Classical spectral approximations, including Jacobi, Laguerre and
Hermite approximations, have gained tremendous success in a variety
of applications, such as Gaussian quadrature rules, root-finding
algorithms and spectral methods for differential and integral
equations (see, e.g.,
\cite{Boyd2000,Davis1984,Shen2011,Trefeth2013}). One of the most
prominent advantages of these approximations is that their
convergence behaviors depend solely on the smoothness properties of
the functions involved and fast convergence can be achieved if the
functions to be approximated are sufficiently smooth. In the past
decade, intensive attention has been devoted to the convergence rate
analysis of Jacobi spectral approximations and their special cases
(i.e., Chebyshev, Legendre, Gegenbauer), and nowadays their
convergence behaviors for analytic and differentiable functions in
weighted and maximum norms are well understood. We refer to
\cite{Liu2019,Trefeth2013,Wang2012,Wang2016,Wang2021,Wang2022,Xiang2020,Xiang2021,Xie2013,Zhao2013}
and references therein for more details. For the latest progress on
error localization of Chebyshev spectral approximations and
pointwise error estimates of Jacobi projections for functions with
singularities, we refer to \cite{Wang2023,Xiang2023}.

Although Laguerre spectral approximations play a pivotal role in
many applications, such as spectral methods for PDEs in unbounded
domain \cite{Shen2000,Shen2009} and the Weeks method for the
numerical inversion of Laplace transform
\cite{Giunta1988,Weeks1966,Weideman1999,Weideman2023}, their
convergence rate analysis has received relatively little attention
compared to their Jacobi counterparts. Even for those existing
results (see, e.g.,
\cite{Elliott1974,Shen2009,Szasz1958,Xiang2012}), however, they are
still far from satisfactory. For example, Elliott and Tuan in
\cite{Elliott1974} studied the estimate of Laguerre coefficients of
analytic functions by using contour integral technique. However, the
contour they used is not the boundary of the convergence domain of
the Laguerre series, which makes it inappropriate for analyzing the
decay rate of the Laguerre coefficients. More recently, Shen and
Wang in \cite{Shen2009} and Xiang in \cite{Xiang2012} studied error
estimates of Laguerre approximations by using the Stirling formula
and the Rodrigues formula for Laguerre polynomials, respectively.
However, the error estimates obtained from these approaches are of
the form $O(n^{-s})$, where $n$ is the number of terms of Laguerre
approximation and $s>0$ is a parameter that depends on the
regularity of the underlying functions. Clearly, such results are
{\it suboptimal} for analytic functions, since from which
exponential convergence cannot be obtained.

In this work we present a thorough convergence rate analysis of
Laguerre spectral approximations for
analytic functions. 
First, we consider the case of Laguerre projection, with the basis
functions being generalized Laguerre polynomials (GLPs) and
generalized Laguerre functions (GLFs). We show that, if the function
to be approximated is analytic inside and on a parabola containing
the positive real axis, then its Laguerre coefficients can be
expressed as a contour integral along the parabola, from which
root-exponential decay of the Laguerre coefficients and
root-exponential convergence of Laguerre projections can be
established. Second, we consider the case of Laguerre interpolation
in Laguerre points and Laguerre-Radau points. By a generalization of
Hermite integral formula for the remainder in polynomial
interpolation, we show that the remainder of Laguerre interpolants
can be expressed as a contour integral along the parabola with which
we then prove that Laguerre interpolants converge at the same
root-exponential rate as their projection counterparts. Finally, we
discuss several applications of Laguerre approximations that are of
practical interest, including Laguerre spectral differentiations,
Gauss-Laguerre quadrature rules, the scaling factors and the Weeks method for the
numerical inversion of Laplace transform. In particular, we obtain
sharp estimates for the convergence rates of Gauss-Laguerre
quadrature rules and the Weeks method for the inversion of Laplace
transform.

We remark that the idea of using contour integral techniques to
analyze the convergence behavior of classical spectral
approximations of analytic functions dates back at least to 1912,
when Bernstein in \cite{Bern1912} proved the exponential convergence
of Chebyshev projections. This idea was adopted in
\cite{Elliott1964,Elliott1974} to study the estimates of Chebyshev,
Jacobi, Laguerre and Hermite coefficients and, more recently, in
\cite{Wang2016,Wang2021,Wang2022} to establish some sharp and
computable upper bounds for the Legendre and Gegenbauer
coefficients. Here we exploit this idea to study optimal convergence
rates of Laguerre projection and interpolation approximations. We
highlight that two key differences between the current study with
\cite{Elliott1974} are that: (i) The contours of the Laguerre
coefficients are different. As will become clear, root-exponential
decay of the Laguerre coefficients can be easily deduced from our
contour, but this is not the case for the contour in \cite{Elliott1974};
(ii) We also study the convergence behavior of Laguerre
interpolation, which was never discussed in \cite{Elliott1974}.

The remainder of this paper is organized as follows. In section
\ref{sec:LagPoly}, we review some useful properties of GLPs and
their weighted Cauchy transform. Sections \ref{sec:LagProj} and
\ref{sec:LagInterp} are devoted to the convergence rate analysis of
Laguerre projection and interpolation approximations, respectively.
In section \ref{sec:Fast} we briefly discuss the convergence
behaviors of Laguerre approximations for entire functions. We
present several important applications of our analysis in section
\ref{sec:Extension} and give some concluding remarks in section
\ref{sec:Conclusion}.

\section{GLPs and their weighted Cauchy transform}\label{sec:LagPoly}
In this section we review some useful properties of GLPs and their
weighted Cauchy transform. For a more comprehensive introduction to
GLPs, we refer to \cite{Ismail2005,Olver2010,Szego1975}.

Let $\mathbb{N}_0:=\mathbb{N}\cup\{0\}$. For each
$k\in\mathbb{N}_0$, the GLP of degree $k$ is defined by
\begin{equation}\label{def:LagPoly}
L_k^{(\alpha)}(x) = \frac{(\alpha+1)_k}{k!} M(-k,\alpha+1;x),
\end{equation}
where $(z)_k$ is the Pochhammer symbol defined by $(z)_{k}
=(z)_{k-1}(z+k-1)$ for $k\in\mathbb{N}$ and $(z)_0=1$, and
$M(a,c;z)$ is the confluent hypergeometric function (also called
Kummer's function) whose definition is given by (see, e.g.,
\cite[Chapter~13]{Olver2010})
\begin{equation}
M(a,c;z) = {}_1 \mathrm{F}_1 \left[\begin{matrix} a & \\
c \end{matrix} \hspace{-.25cm} ;  x \right] = \sum_{k=0}^{\infty}
\frac{(a)_k}{(c)_k} \frac{z^k}{k!}. \notag
\end{equation}
In the particular case of $\alpha=0$, we shall drop the superscript
of GLP and write $L_k^{(0)}(x)=L_k(x)$ for simplicity. Let
$\omega_{\alpha}(x)=x^{\alpha}e^{-x}$, where $\alpha>-1$, be the
generalized Laguerre weight function and let
$\mathbb{R}_{+}:=[0,\infty)$. It is well known that the GLPs are
orthogonal with respect to the inner product
\begin{equation}\label{def:InnerProd}
\langle f,g \rangle_{\omega_{\alpha}} = \int_{\mathbb{R}_{+}} f(x)
g(x) \omega_{\alpha}(x) \mathrm{d}x.  
\end{equation}
More precisely, we have
\begin{equation}
\langle L_{k}^{(\alpha)},L_{j}^{(\alpha)} \rangle_{\omega_{\alpha}}
= \gamma_{k}^{(\alpha)}\delta_{k,j}, \quad \gamma_{k}^{(\alpha)}=
\frac{\Gamma(k+\alpha+1)}{k!}, \nonumber
\end{equation}
where $\delta_{k,j}$ is the Kronecker delta and $\Gamma(z)$ is the
gamma function. Moreover, using the asymptotic expansion of the
ratio of gamma functions \cite[Equation~(5.11.13)]{Olver2010}, we
find for large $k$ that
\begin{equation}\label{eq:gammaAsym}
\gamma_{k}^{(\alpha)} = k^{\alpha} \left[1 +
\frac{\alpha(\alpha+1)}{2k} + O(k^{-2}) \right].
\end{equation}
For $x\in\mathbb{R}_{+}$, the GLPs satisfy the following inequality
\begin{equation}\label{eq:LagFunBound}
e^{-x/2} |L_n^{(\alpha)}(x)| \leq \left\{
\begin{array}{ll}
{\displaystyle \kappa_n^{\alpha}},    &  \alpha\geq0,  \\[10pt]
{\displaystyle 2-\kappa_n^{\alpha}},  &  \alpha\in(-1,0),
\end{array}
\right. \quad
\kappa_n^{\alpha}=\frac{\gamma_{n}^{(\alpha)}}{\Gamma(\alpha+1)},
\end{equation}
which will be useful in proving error estimates of Laguerre
approximations using GLFs.

Let $z\in\mathbb{C}\setminus\mathbb{R}_{+}$, the weighted Cauchy
transform of GLPs is defined by
\begin{equation}\label{def:Phi}
\Phi_n^{(\alpha)}(z) = \frac{1}{2\pi\mathrm{i}}\int_{\mathbb{R}_{+}}
\frac{\omega_{\alpha}(x)L_n^{(\alpha)}(x)}{z-x} \mathrm{d}x,
\end{equation}
where $\mathrm{i}$ is the imaginary unit.
Below we summarize several important properties satisfied by
$\Phi_n^{(\alpha)}(z)$.
\begin{lemma}
The following statements are true.
\begin{itemize}
\item[\rm (i)] $\Phi_n^{(\alpha)}(z)$ is analytic in the
whole complex plane cut along $\mathbb{R}_{+}$.

\item[\rm (ii)] For $x\in(0,\infty)$, it holds that
\begin{align}
\lim_{\varepsilon\rightarrow0^{+}} \left[
\Phi_n^{(\alpha)}(x-\mathrm{i}\varepsilon) +
\Phi_n^{(\alpha)}(x+\mathrm{i}\varepsilon) \right] &= \frac{1}{\pi
\mathrm{i}}\dashint_{\mathbb{R}_{+}} \frac{\omega_{\alpha}(y)
L_n^{(\alpha)}(y)}{x-y} \mathrm{d}y, \label{eq:PlemeljI}  \\
\lim_{\varepsilon\rightarrow0^{+}} \left[
\Phi_n^{(\alpha)}(x-\mathrm{i}\varepsilon) -
\Phi_n^{(\alpha)}(x+\mathrm{i}\varepsilon) \right] &=
\omega_{\alpha}(x)L_n^{(\alpha)}(x), \label{eq:PlemeljII}
\end{align}
where the bar indicates the Cauchy principal value.

\item[\rm (iii)] As $z\rightarrow0$, it holds that
\begin{equation}\label{eq:PhiAsym}
\Phi_n^{(\alpha)}(z) = \left\{
\begin{array}{ll}
{\displaystyle O(1)},           &  \alpha>0,  \\[8pt]
{\displaystyle O(\log z)},      &  \alpha=0,  \\[8pt]
{\displaystyle O(z^{\alpha})},  &  \alpha<0.
\end{array}
\right.
\end{equation}
Moreover, $\Phi_n^{(\alpha)}(z)=O(z^{-n-1})$ as
$z\rightarrow\infty$.

\item[\rm (iv)] It holds that
\begin{align}\label{eq:PhiU}
\Phi_n^{(\alpha)}(z) &= \frac{\mathrm{i}}{2\pi}\Gamma(n+\alpha+1)
U(n+1,1-\alpha; -z),
\end{align}
where $U(a,c;z)$ is the confluent hypergeometric function or
Kummer's function of the second kind (see, e.g.,
\cite[Chapter~13]{Olver2010}).

\item[\rm (v)] As $n\rightarrow\infty$, it holds that
\begin{align}\label{eq:PhiAsym}
\Phi_n^{(\alpha)}(z) &= \frac{\mathrm{i}\gamma_n^{(\alpha)}
e^{-2\sqrt{-(n+1)z}}}{2\sqrt{\pi}}\left[
\frac{e^{-z/2}(-z)^{\alpha/2+1/4}}{(n+1)^{\alpha/2+1/4}} +
O(n^{-\alpha/2-3/4}) \right].
\end{align}
\end{itemize}
\end{lemma}
\begin{proof}
The first statement is trivial and the second statement follows from
the well-known Plemelj formula\footnote{It also called the
Sokhotskyi-Plemelj formula in some literature.}
\cite[Chapter~2]{Muskhl1972}. As for the third statement, the
asymptotic behavior as $z\rightarrow0$ follows from the
matrix-valued Riemann-Hilbert problem with respect to the case of
GLPs \cite[Chapter~22]{Ismail2005} and the asymptotic behavior as
$z\rightarrow\infty$ follows from the orthogonality of the GLPs. For
the fourth statement, it was given in \cite[Table~1]{Elliott1974}
but without proof. Below we provide a concise proof for the purpose
of self-contained. Combining the connection formula between GLPs and
Whittaker functions \cite[Equation~(18.11.2)]{Olver2010} and
\cite[Equation~(7.627.8)]{Grad2007}, we obtain that
\begin{align}\label{eq:PhiU1}
\Phi_n^{(\alpha)}(z) 
&= -\frac{\Gamma(n+\alpha+1)}{2\pi\mathrm{i}} \left[
\frac{\Gamma(\alpha)}{\Gamma(n+\alpha+1)} M(n+1,1-\alpha,-z) \right. \nonumber \\
&~~~~~~~~~ \left. + \frac{\Gamma(-\alpha)
(-z)^{\alpha}}{\Gamma(n+1)} M(n+\alpha+1,1+\alpha,-z) \right].
\end{align}
Let $\Upsilon$ denote the term inside the square brackets. In the
case where $\alpha\notin\mathbb{N}_{0}$, using
\cite[Equation~(13.2.42)]{Olver2010}, we find that $\Upsilon =
U(n+1,1-\alpha;-z)$. Substituting this into \eqref{eq:PhiU1} gives
\eqref{eq:PhiU}. This proves the case $\alpha\notin\mathbb{N}_{0}$.
On the other hand, using the reflection formula
\cite[Equation~(5.5.3)]{Olver2010}, we obtain that
\begin{align}\label{eq:PhiU2}
\Upsilon &= \frac{\pi}{\sin(\alpha\pi)} \left[
\frac{\mathrm{\mathbf{M}}(n+1,1-\alpha; -z)}{\Gamma(n+\alpha+1)} -
\frac{(-z)^{\alpha}}{\Gamma(n+1)}
\mathrm{\mathbf{M}}(n+\alpha+1,\alpha+1;-z) \right].
\end{align}
where $\mathrm{\mathbf{M}}(a,c;z)=M(a,c;z)/\Gamma(c)$, which is an
entire function of $a,c,z$ (see, e.g.,
\cite[Equation~(13.2.4)]{Olver2010}). Furthermore, using
\cite[Equation~(13.2.5)]{Olver2010} it is easily verified that the
term in square brackets on the right-hand side of \eqref{eq:PhiU2}
tends to zero as $\alpha$ tends to a nonnegative integer, and hence
$\Upsilon = U(n+1,1-\alpha;-z)$ still holds in the case of
$\alpha\in\mathbb{N}_0$ if we take the limit by l'H\^{o}pital's
rule. This proves the fourth statement. As for the last statement,
combining \cite[Equation~(13.8.8)]{Olver2010} and
\cite[Equation~(10.40.2)]{Olver2010} we get
\begin{equation}\label{def:KummerUAsy}
U(a,c;z) = \frac{\sqrt{\pi}}{\Gamma(a)}
\left(\frac{z}{a}\right)^{(1-c)/2}
\frac{\exp(z/2-2\sqrt{az})}{(az)^{1/4}} \left(1+O(a^{-1/2})\right),
\end{equation}
as $a\rightarrow\infty$. The desired result then follows from
applying the above asymptotic estimate to \eqref{eq:PhiU}. This ends
the proof.
\end{proof}

\section{Convergence rate analysis of Laguerre projections}\label{sec:LagProj}
In this section, we conduct a thorough convergence rate analysis of
Laguerre projections. Throughout the paper, we always take the
positive orientation for a contour integral is the counterclockwise
orientation. Moreover, we denote by $\mathbb{P}_n$ the space of
polynomials of degree up to $n$ and denote by $\mathcal{K}$ a
generic positive constant 
and may take different values at different places.

Before starting our analysis, we introduce a parabola in the complex
plane
\begin{equation}\label{def:Parab}
P_{\rho} := \left\{ z\in\mathbb{C}: ~\Re(\sqrt{-z})=\rho \right\},
\end{equation}
where $\rho>0$ and the branch cut is taken along the positive real
axis and we choose the branch such that $\sqrt{-z}$ is taken real
and positive if $z<0$. Moreover, we denote by $D_{\rho}$ the
interior domain of $P_{\rho}$, i.e.,
\begin{equation}\label{def:D}
D_{\rho}=\{z\in\mathbb{C}: ~\Re(\sqrt{-z})<\rho\}.
\end{equation}
From \eqref{def:Parab} one can easily check that $P_{\rho}$ is a
parabola opening to the right with vertex at $z=-\rho^2$ and focus
at the origin, and it degenerates into the positive real axis as
$\rho\rightarrow0$. A key motivation for introducing $P_{\rho}$ is
that the boundary of the convergence domain of Laguerre
approximations of analytic functions can be characterized by
$P_{\rho}$ for some $\rho>0$ (see \cite[Chapter~9]{Szego1975} and
\cite{Szasz1958}).

\begin{figure}[ht]
\includegraphics[width=0.45\textwidth]{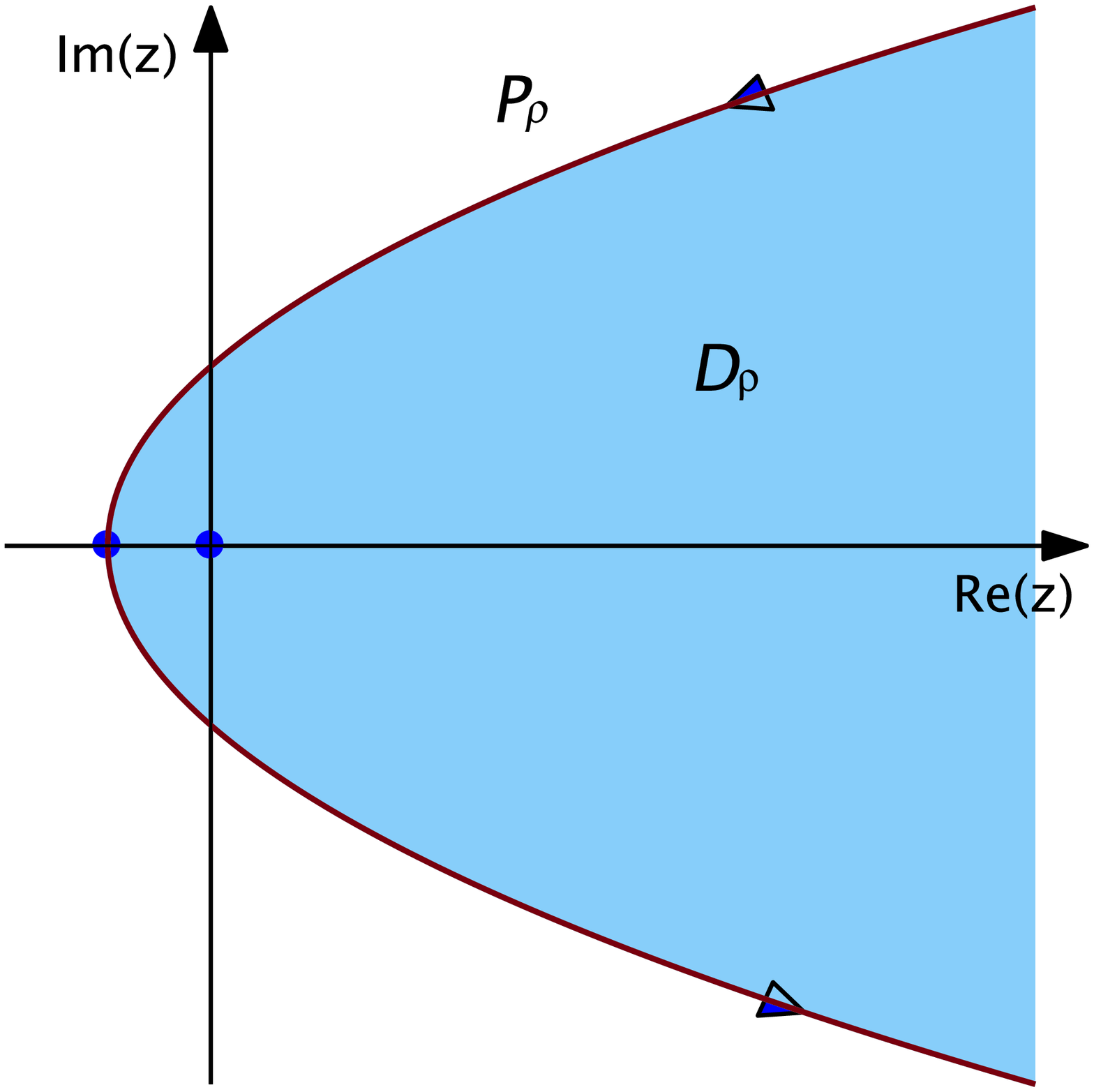}\qquad
\includegraphics[width=0.45\textwidth]{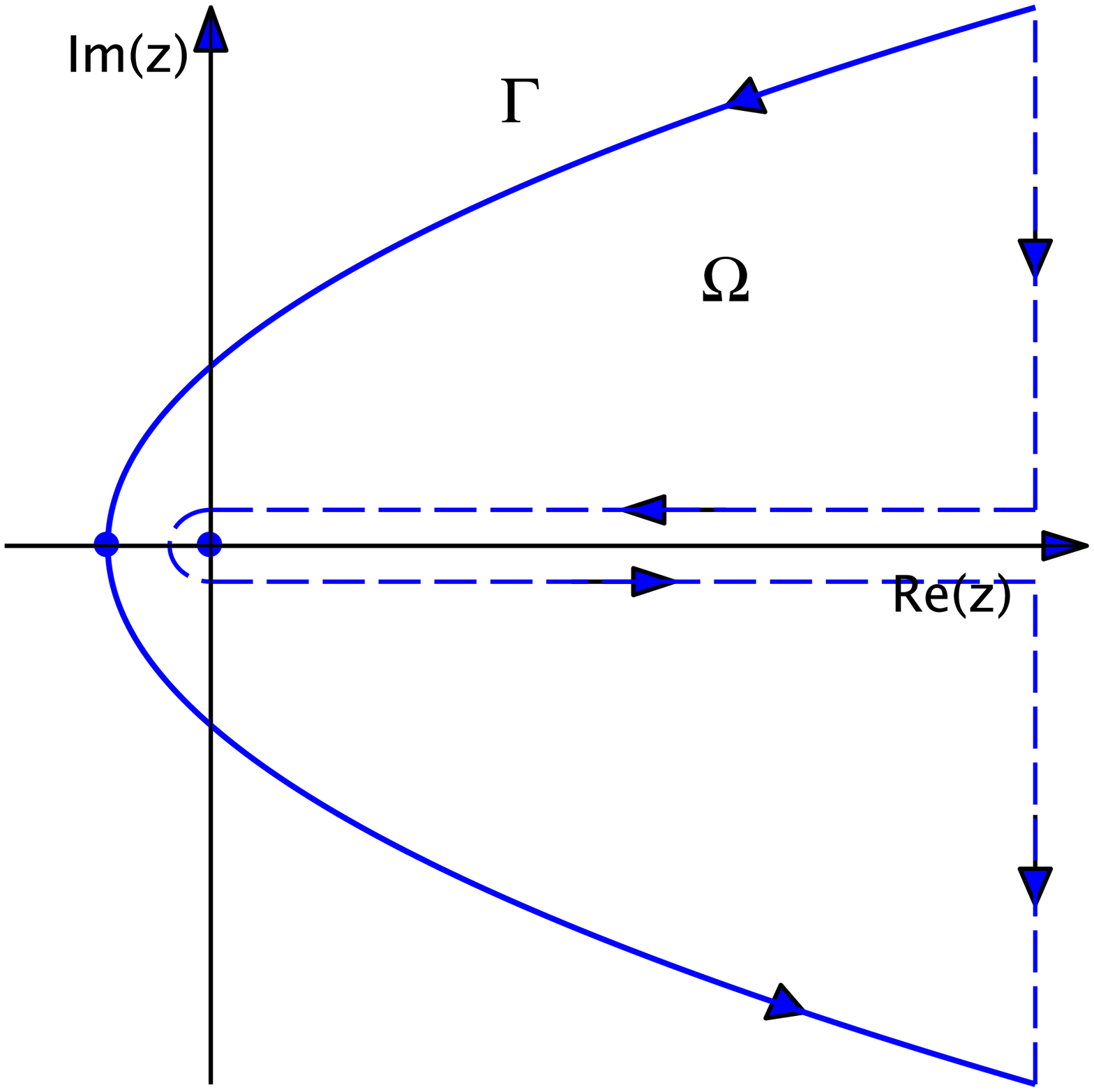}
\caption{Left: The parabola $P_{\rho}$ and the region $D_{\rho}$.
The two points on the horizontal axis (from left to right) are the
vertex and focus of $P_{\rho}$, respectively. Right: The dashed
lines indicates the five contours used in the proof of Theorem
\ref{thm:LagParab}.}\label{fig:Parab}
\end{figure}

\subsection{Projections using GLPs}
Let $L_{\omega_{\alpha}}^2(\mathbb{R}_{+})$ denote the Hilbert space
consisting of all square integrable with respect to the weight
function $\omega_{\alpha}(x)$ on $\mathbb{R}_{+}$ and let
$\Pi_n^\mathrm{P}$ denote the orthogonal projection operator from
${L}_{\omega_{\alpha}}^2(\mathbb{R}_{+})$ upon the polynomial space
$\mathbb{P}_n$. Then, for any
$f\in{L}_{\omega_{\alpha}}^2(\mathbb{R}_{+})$ we have
\begin{equation}\label{def:LagProj}
(\Pi_n^\mathrm{P}f)(x) = \sum_{k=0}^{n} a_k^{(\alpha)}
L_k^{(\alpha)}(x), \quad a_k^{(\alpha)} = \frac{\langle
f,L_k^{(\alpha)} \rangle_{\omega_{\alpha}}}{\gamma_{k}^{(\alpha)}}.
\end{equation}
We now state the first main result of this work.
\begin{theorem}\label{thm:LagParab}
Assume that $f$ is analytic within and on the parabola $P_{\rho}$
for some $\rho>0$ and $|f(z)|\leq\mathcal{K}|z|^{\beta}$ for some
$\beta\in\mathbb{R}$ as $z\rightarrow\infty$ in $D_{\rho}$. Then,
for each $k\geq\max\{\lfloor\beta-1/2\rfloor+1,0\}$ where $\lfloor
\cdot\rfloor$ denotes the integer part,
\begin{equation}\label{eq:LagCont}
a_k^{(\alpha)} = \frac{1}{\gamma_{k}^{(\alpha)}} \int_{P_{\rho}}
\Phi_k^{(\alpha)}(z) f(z) \mathrm{d}z.
\end{equation}
\end{theorem}
\begin{proof}
Let $\eta>0$ be arbitrary and let $\Omega$ denote the region
enclosed by $P_{\rho}$ and the vertical line $\Re(z)=\eta$.
Moreover, let $\Gamma$ denote the overlapping boundaries of $\Omega$
and $D_{\rho}$. We consider the integral of
$\Phi_k^{(\alpha)}(z)f(z)$ along the contour $\Gamma$. Since
$\Phi_k^{(\alpha)}(z)f(z)$ is analytic in $\Omega$ except along the
interval $[0,\eta]$, from Cauchy's Theorem we can deform the contour
$\Gamma$ into the following contours (see the right panel of Figure
\ref{fig:Parab}):
\begin{itemize}
\item[(i)] $\Gamma_V^{1}$: the vertical line segment from the upper
intersection point between $P_{\rho}$ and the vertical line
$\Re(z)=\eta$ to the point $z=\eta+\mathrm{i}\varepsilon$, where
$\varepsilon>0$ is small;

\item[(ii)] $\Gamma_H^{1}$: the horizontal line segment from
$z=\eta+\mathrm{i}\varepsilon$ to $z=\mathrm{i}\varepsilon$;

\item[(iii)] $\Gamma_\varepsilon$: the semicircle centered at the
origin with radius $\varepsilon$ in the left half plane;

\item[(iv)] $\Gamma_H^{2}$: the horizontal line segment from
$z=-\mathrm{i}\varepsilon$ to $z=\eta-\mathrm{i}\varepsilon$;

\item[(v)] $\Gamma_V^{2}$: the vertical line segment from
$z=\eta-\mathrm{i}\varepsilon$ to the lower intersection point
between $P_{\rho}$ and the vertical line $\Re(z)=\eta$.
\end{itemize}
Thus
\begin{equation}
\int_{\Gamma} \Phi_k^{(\alpha)}(z) f(z) \mathrm{d}z = \left(
\int_{\Gamma_V^{1}} + \int_{\Gamma_H^{1}} +
\int_{\Gamma_\varepsilon} + \int_{\Gamma_H^{2}} +
\int_{\Gamma_V^{2}} \right) \Phi_k^{(\alpha)}(z) f(z) \mathrm{d}z.
\nonumber
\end{equation}
In the following we shall consider the limits of these five contour
integrals on the right-hand side as $\varepsilon\rightarrow0^{+}$
and $\eta\rightarrow\infty$. We start with the third contour
integral. By \eqref{eq:PhiAsym} we immediately get
\begin{align}
\left|\int_{\Gamma_\varepsilon} \Phi_k^{(\alpha)}(z)f(z) \mathrm{d}z
\right| &= \left\{
\begin{array}{ll}
{\displaystyle O(\varepsilon)},                  &  \alpha>0,  \\[8pt]
{\displaystyle O(\varepsilon\log\varepsilon)},   &  \alpha=0,  \\[8pt]
{\displaystyle O(\varepsilon^{\alpha+1})},       &  \alpha<0,
\end{array}
\right. \nonumber
\end{align}
and therefore the third contour integral vanishes as
$\varepsilon\rightarrow0^{+}$. For the second and fourth contour
integrals, we find that
\begin{align}
\lim_{\substack{\varepsilon\rightarrow0^{+} \\
\eta\rightarrow\infty}} \left(\int_{\Gamma_H^{1}} +
\int_{\Gamma_H^{2}} \right) \Phi_k^{(\alpha)}(z)f(z) \mathrm{d}z &=
\lim_{\varepsilon\rightarrow0^{+}} \int_{\mathbb{R}_{+}} \left[
\Phi_k^{(\alpha)}(x-\mathrm{i}\varepsilon) -
\Phi_k^{(\alpha)}(x+\mathrm{i}\varepsilon) \right] f(x) \mathrm{d}x \notag \\
&= \int_{\mathbb{R}_{+}}\omega_{\alpha}(x)L_k^{(\alpha)}(x) f(x) \mathrm{d}x,  \nonumber \\
&= a_k^{(\alpha)} \gamma_k^{(\alpha)}, \nonumber
\end{align}
where we have used the Plemelj formula \eqref{eq:PlemeljII} in the
second step. For the first and last contour integrals, we obtain
after some calculation
\begin{align}
\lim_{\substack{\varepsilon\rightarrow0^{+} \\
\eta\rightarrow\infty}}\left( \int_{\Gamma_V^{1}} +
\int_{\Gamma_V^{2}} \right) \Phi_k^{(\alpha)}(z)f(z) \mathrm{d}z &=
(-\mathrm{i}) \lim_{\eta\rightarrow\infty}
\int_{-2\rho\sqrt{\eta+\rho^2}}^{2\rho\sqrt{\eta+\rho^2}}
\Phi_k^{(\alpha)}(\eta+\mathrm{i}y) f(\eta+\mathrm{i}y) \mathrm{d}y.
\nonumber
\end{align}
Recalling that $|\Phi_k^{(\alpha)}(z)|=O(|z|^{-k-1})$ and
$|f(z)|=O(|z|^{\beta})$ as $|z|\rightarrow\infty$ in $D_{\rho}$, one
can easily check that the integral on the right-hand side of the
above equation vanishes as $\eta\rightarrow\infty$ for
$k\geq\max\{\lfloor\beta-1/2\rfloor+1,0\}$. Combining these results
and noting that $\Gamma\rightarrow{P}_{\rho}$ as
$\eta\rightarrow\infty$, the desired result follows immediately.
This completes the proof.
\end{proof}

\begin{remark}\label{rk:contour}
Elliott and Tuan in \cite[Theorem~1]{Elliott1974} also derived a
contour integral representation for the Laguerre coefficients
$\{a_k^{(\alpha)}\}_{k=0}^{\infty}$ and the contour was defined as a
combination of two connected contours with one lying in the upper
half plane which connects $z=\infty$ to $z=0$ and the other lying in
the lower half planes which connects $z=0$ to $z=\infty$, and the
slopes of both connected contours tend to zero as $z\rightarrow0$
and $z\rightarrow\infty$, respectively. Clearly, this contour is
different from our contour in Theorem \ref{thm:LagParab}. Although Elliott and Tuan
derived the exact formula for the Laguerre coefficients of the function $f(x)=1/(x+\lambda)$ with $\lambda>0$,
the analysis of the decay behavior of the Laguerre coefficients with this contour is highly nontrivial. In contrast, a root-exponential decay of the Laguerre coefficients can be deduced from our contour immediately, as will be shown below.
\end{remark}

\begin{remark}\label{rk:contour}
Note that \eqref{eq:LagCont} holds for all $k\geq0$ whenever
$\beta<1/2$.
\end{remark}

We now state the second main result of this work, which establishes
a root-exponential convergence of Laguerre projections. Hereafter,
we will denote by $\mathrm{d}s$ the differential of arc length and
denote by $\|\cdot\|_{\omega_{\alpha}}$ the norm induced by the
inner product \eqref{def:InnerProd}.
\begin{theorem}\label{thm:LagPolyBound}
Under the assumptions of Theorem \ref{thm:LagParab} and let
$V_{\alpha}$ be defined by
\[
V_{\alpha} =
\int_{P_{\rho}}|(-z)^{\alpha/2+1/4}e^{-z/2}f(z)|\mathrm{d}s.
\]
Then, for each $n\geq\max\{\lfloor\beta-1/2\rfloor+1,0\}$, the
following statements hold.
\begin{itemize}
\item[\rm(i)] The Laguerre coefficient satisfies
\begin{equation}\label{eq:LagCoeffBound}
|a_n^{(\alpha)}| \leq \frac{\mathcal{K}
e^{-2\rho\sqrt{n+1}}}{(n+1)^{\alpha/2+1/4}},
\end{equation}
where $\mathcal{K}\cong V_{\alpha}/(2\sqrt{\pi})$ for $n\gg1$. 

\item[\rm(ii)] The error of the Laguerre projection satisfies
\begin{align}\label{eq:LagProjBound}
\left\|f - \Pi_n^\mathrm{P}f\right\|_{\omega_{\alpha}} \leq
\mathcal{K}e^{-2\rho\sqrt{n+1}},
\end{align}
where $\mathcal{K}\cong V_{\alpha}/(2\sqrt{2\pi\rho})$ for $n\gg1$. 
\end{itemize}
\end{theorem}
\begin{proof}
We first consider the proof of \eqref{eq:LagCoeffBound}. From
Theorem \ref{thm:LagParab} it follows that
\begin{equation}
|a_n^{(\alpha)}| \leq \frac{1}{\gamma_{n}^{(\alpha)}}
\int_{P_{\rho}} |\Phi_n^{(\alpha)}(z)f(z)| \mathrm{d}s.  \nonumber
\end{equation}
On the other hand, from \eqref{eq:PhiAsym} we obtain for
$z\in{P}_{\rho}$ that
\begin{align}
|\Phi_n^{(\alpha)}(z)| &=
\frac{\gamma_n^{(\alpha)}e^{-2\rho\sqrt{n+1}}}{2\sqrt{\pi}} \left[
\frac{|e^{-z/2}(-z)^{\alpha/2+1/4}|}{(n+1)^{\alpha/2+1/4}} +
O(n^{-\alpha/2-3/4}) \right]. \nonumber
\end{align}
Hence the inequality \eqref{eq:LagCoeffBound} follows by combining
the above two results. We next consider the proof of
\eqref{eq:LagProjBound}. By \eqref{eq:LagCoeffBound} and the
estimate of $\gamma_n^{(\alpha)}$ in \eqref{eq:gammaAsym}, we have
\begin{align}
\|f-\Pi_n^\mathrm{P}f\|_{\omega_{\alpha}}^2 = \sum_{k=n+1}^{\infty}
(a_k^{(\alpha)})^2 \gamma_k^{(\alpha)} &\leq \mathcal{K}^2
\sum_{k=n+1}^{\infty} \frac{e^{-4\rho\sqrt{k+1}}}{\sqrt{k+1}},
\nonumber
\end{align}
where $\mathcal{K}\cong V_{\alpha}/(2\sqrt{\pi})$. 
Furthermore, the last sum can be bounded by
\begin{align}
\sum_{k=n+1}^{\infty} \frac{e^{-4\rho\sqrt{k+1}}}{\sqrt{k+1}} \leq
\int_{n+1}^{\infty} \frac{e^{-4\rho\sqrt{x}}}{\sqrt{x}} \mathrm{d}x
= \frac{e^{-4\rho\sqrt{n+1}}}{2\rho}. \nonumber
\end{align}
Combining the above two inequalities gives \eqref{eq:LagProjBound}.
This ends the proof.
\end{proof}

\begin{remark}
When $\max_{z\in P_{\rho}}|f(z)|$ is bounded, then the quantity $V_{\alpha}$ can be further bounded by
\begin{align}
V_{\alpha} &\leq \max_{z\in P_{\rho}}|f(z)| \int_{P_{\rho}} |(-z)^{\alpha/2+1/4}e^{-z/2}|\mathrm{d}s \nonumber \\
&= 2 \max_{z\in P_{\rho}}|f(z)| \int_{-\infty}^{\infty} (t^2+\rho^2)^{\alpha/2+3/4} e^{-(t^2-\rho^2)/2} \mathrm{d}t \nonumber \\
&= 2 \max_{z\in P_{\rho}}|f(z)| \sqrt{\pi} e^{\rho^2/2} \rho^{\alpha+5/2} U(1/2,\alpha/2+9/4;\rho^2/2), \nonumber
\end{align}
where we have parametrized $P_{\rho}$ by $z(t)=t^2-\rho^2-2t\rho\mathrm{i}$ with $t\in(-\infty,\infty)$ in the second step.
\end{remark}

\begin{figure}[ht]
\centering
\includegraphics[width=0.49\textwidth,height=0.42\textwidth]{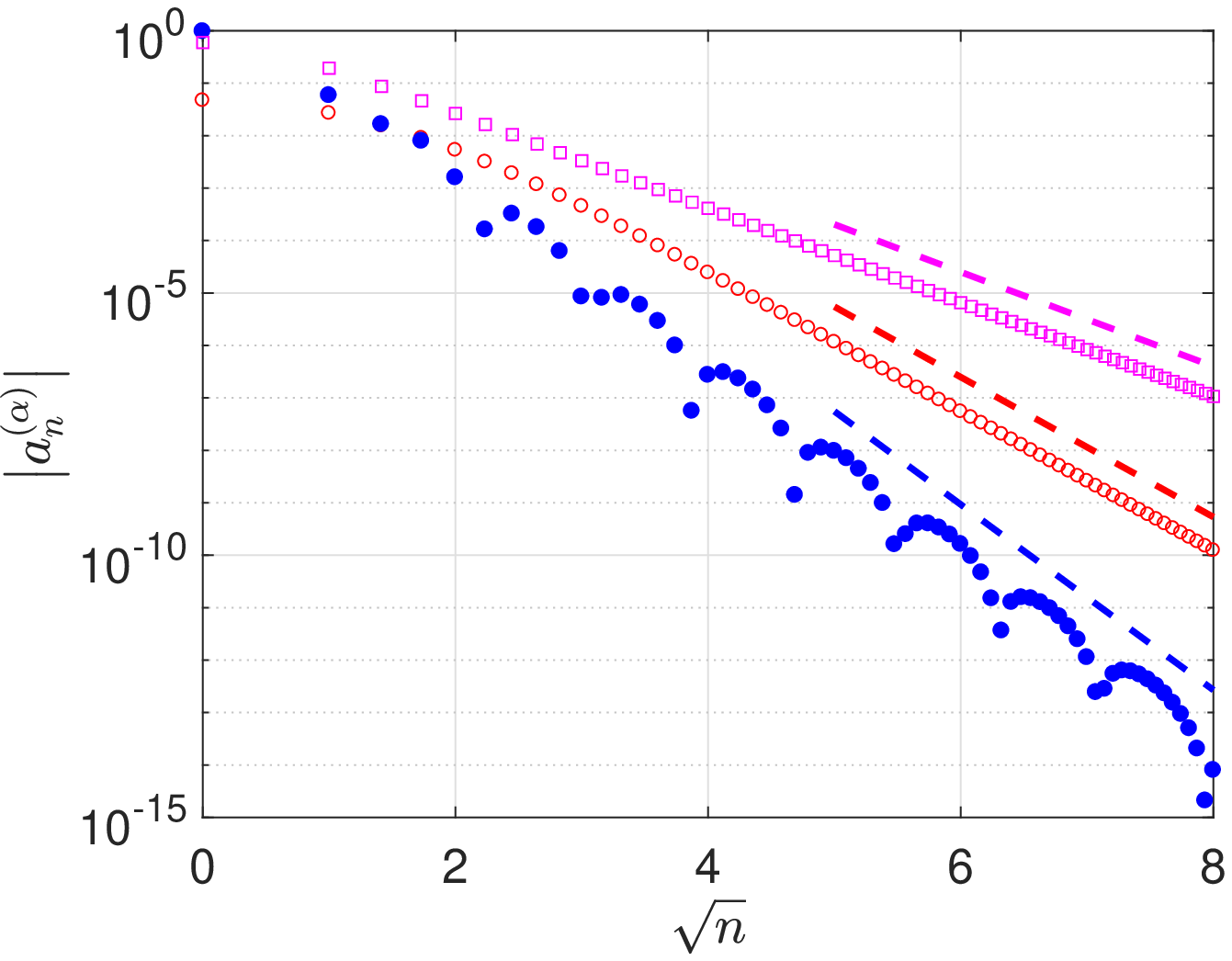}~
\includegraphics[width=0.49\textwidth,height=0.42\textwidth]{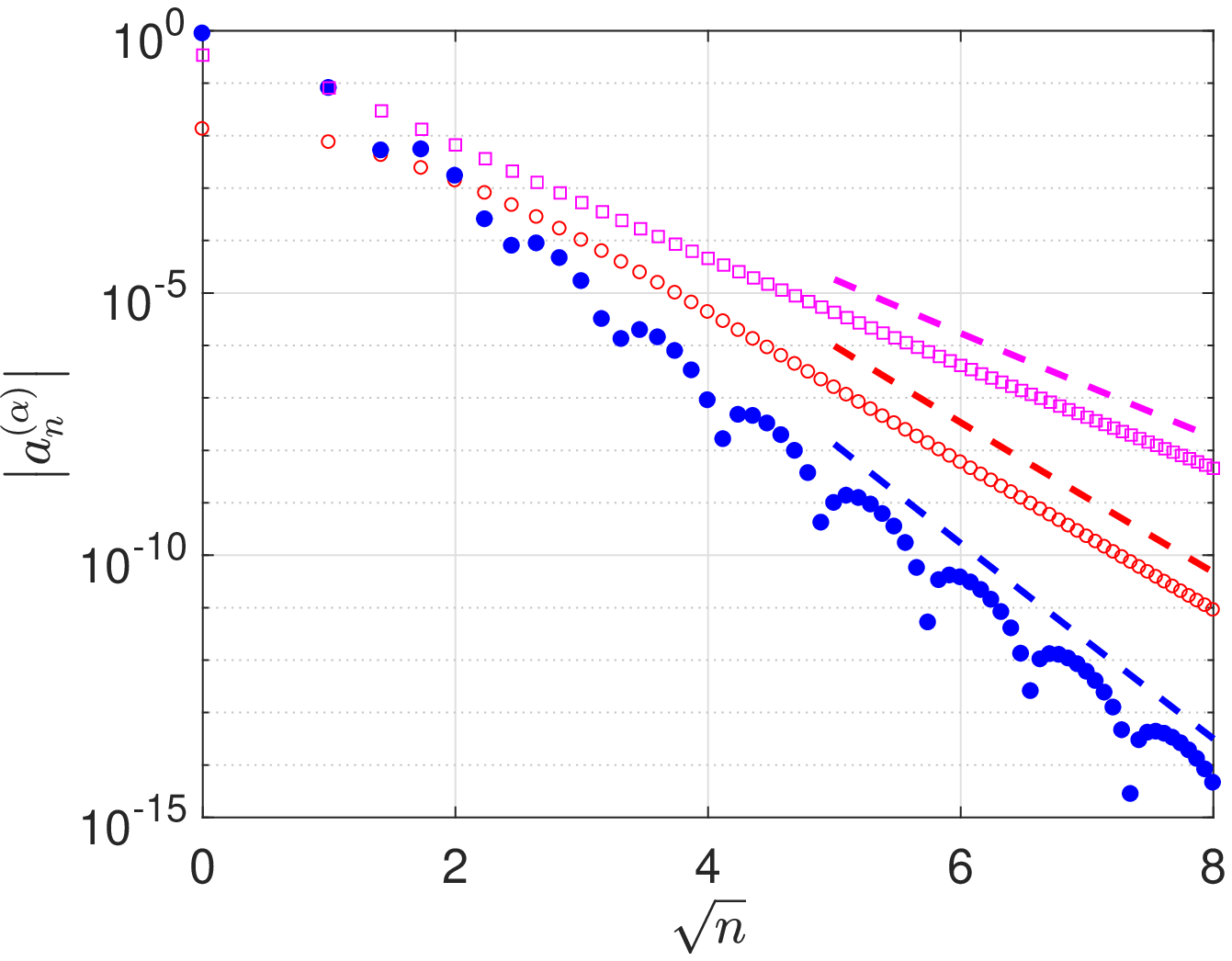}
\caption{The magnitudes of the Laguerre coefficients of $f_1(x)=1/(x+1)$
(squares), $f_2(x)=e^{-x}/(4x+9)$ (circles) and $f_3(x)=\mathrm{sech}(\pi{x}/16)$ (dots). Here $\alpha=0$ (left)
and $\alpha=3/2$ (right) and the dashed lines show the corresponding
predicted rates $O(n^{-\alpha/2-1/4}e^{-2\rho\sqrt{n}})$.}
\label{fig:Lagcoeff}
\end{figure}

\begin{remark}\label{rk:Pararho}
If $f(z)$ has branch point or pole singularities $\{z_k\}_{k=1}^{m}$
in the complex plane, with each
$z_k\in\mathbb{C}\setminus\mathbb{R}_{+}$, then it is easily
verified that $f$ is analytic within the parabola $P_{\rho}$, where
$\rho$ is defined by
\begin{equation}\label{eq:Pararho}
\rho = \min_{k=1,\ldots,m} \left\{ \sqrt{\frac{|z_k| - \Re(z_k)}{2}}
\right\} -\varepsilon,
\end{equation}
where $\varepsilon$ is arbitrarily close to zero. Moreover, from Theorem \ref{thm:LagPolyBound} we conclude that
$|a_n^{(\alpha)}|=O(n^{-\alpha/2-1/4}e^{-2\rho\sqrt{n}})$ as
$n\rightarrow\infty$ and thus the Laguerre coefficients decay at a
root-exponential rate; see Figure \ref{fig:Lagcoeff} for an
illustration.
\end{remark}

\begin{remark}
If $f(z)$ is analytic inside the parabolic region $D_{\rho}$ with
$\rho>0$ and if there exists a positive number
$\mathcal{B}(\alpha,\rho)$ such that
\[
|f(z)|\leq \mathcal{B}(\alpha,\rho)\exp\left(\frac{\Re(z)}{2} -
\sqrt{|\Re(z)|\left(\rho^2 - \frac{|z|-\Re(z)}{2}\right)}\right),
\quad  z\in D_{\rho}.
\]
Then, for $\alpha\in(-1,1/2)$, it has been shown in
\cite[Lemma~3.4]{Szasz1958} that the Laguerre coefficients satisfy
$|a_n^{(\alpha)}|\leq \mathcal{B}(\rho)e^{-2\rho\sqrt{n}}$ for some
$\mathcal{B}(\rho)>0$. Comparing this with \eqref{eq:LagCoeffBound},
we can see that \eqref{eq:LagCoeffBound} is better since it holds
for all $\alpha>-1$ and it is more precise.
\end{remark}

Figure \ref{fig:Lagcoeff} illustrates the magnitudes of the Laguerre
coefficients $a_n^{(\alpha)}$ on a log scale with $\sqrt{n}$ on the
horizontal axis of the test functions $f_1(x)=1/(x+1)$,
$f_2(x)=e^{-x}/(4x+9)$ and $f_3(x)=\mathrm{sech}(\pi x/16)$. It is
easily seen that $f_1$ and $f_2$ have simple poles at $z=-1$ and
$z=-9/4$, respectively, and $f_3$ has poles at
$z=\pm8(2k+1)\mathrm{i}$ with $k=0,1,\ldots$. Moreover, it is not
difficult to verify that they all satisfy the assumptions of Theorem
\ref{thm:LagParab} and therefore, by Remark \ref{rk:Pararho}, the
predicted decay rates of their Laguerre coefficients are
$O(n^{-\alpha/2-1/4}e^{-2\rho\sqrt{n}})$ with $\rho=1,3/2,2$,
respectively. From Figure \ref{fig:Lagcoeff} we can see that the
actual decay rates of the Laguerre coefficients are in good
agreement with the predicted ones.

\subsection{Projections using GLFs}
In some contexts, the use of GLPs may be inconvenient since they are
unbounded on the positive real axis. To remedy this, it is
preferable to use the GLFs which are defined by
\begin{equation}\label{def:LagFun}
\widehat{L}_k^{(\alpha)}(x) = e^{-x/2} {L}_k^{(\alpha)}(x), \quad
k=0,1,\ldots.
\end{equation}
We introduce the inner product and the associated norm
\begin{equation}
\langle f,g \rangle_{\varpi_{\alpha}} = \int_{\mathbb{R}_{+}} f(x)
g(x)\varpi_{\alpha}(x) \mathrm{d}x, \quad \|f\|_{\varpi_{\alpha}} =
\sqrt{\langle f,f \rangle_{\varpi_{\alpha}}}, \nonumber
\end{equation}
where $\varpi_{\alpha}(x)=x^{\alpha}$ with $\alpha>-1$. Denote by
$L_{\varpi_{\alpha}}^2(\mathbb{R}_{+})$ the space of functions that
are square integrable on $\mathbb{R}_{+}$ with respect to the weight
function $\varpi_{\alpha}(x)$. It is easy to check that the GLFs are
orthogonal with respect to the above inner product and
\[
\langle\widehat{L}_k^{(\alpha)},\widehat{L}_j^{(\alpha)}
\rangle_{\varpi_{\alpha}} = \gamma_k^{(\alpha)}\delta_{k,j}.
\]
Let $\mathbb{Q}_n=\mathrm{span}\{e^{-x/2}u(x),~u\in\mathbb{P}_n\}$
and let $\Pi_n^{\mathrm{F}}$ denote the orthogonal projection
operator from $L_{\varpi_{\alpha}}^2(\mathbb{R}_{+})$ upon the
function space $\mathbb{Q}_n$. Then, for any
$f\in{L}_{\varpi_{\alpha}}^2(\mathbb{R}_{+})$ we have
\begin{equation}\label{def:LagProj2}
(\Pi_n^{\mathrm{F}}f)(x) = \sum_{k=0}^{n} b_k^{(\alpha)}
\widehat{L}_k^{(\alpha)}(x),  \quad  b_k^{(\alpha)} = \frac{\langle
f,\widehat{L}_k^{(\alpha)}\rangle_{\varpi_{\alpha}}}{\gamma_{k}^{(\alpha)}}.
\end{equation}
Below we shall establish a sharp estimate on the decay rate of the
Laguerre coefficients $\{b_k^{(\alpha)}\}_{k=0}^{\infty}$ and error
estimates of $\Pi_n^{\mathrm{F}}f$ in the weighted and maximum
norms.
\begin{theorem}\label{thm:LagFunBound}
Let $f$ be analytic within and on the parabola $P_{\rho}$ for some
$\rho>0$ and let $|f(z)|\leq\mathcal{K}|z|^{\beta}e^{-\Re(z)/2}$ for
some $\beta\in\mathbb{R}$ as $|z|\rightarrow\infty$ in $D_{\rho}$.
Moreover, let $\widehat{V}_{\alpha}$ be defined by
\[
\widehat{V}_{\alpha} = \int_{P_{\rho}}|(-z)^{\alpha/2+1/4}f(z)|\mathrm{d}s.
\]
Then, for each $n\geq\max\{\lfloor\beta-1/2\rfloor+1,0\}$, the
following statements hold.
\begin{itemize}
\item[\rm(i)] The Laguerre coefficient satisfies
\begin{equation}\label{eq:akF}
|b_n^{(\alpha)}| \leq
\frac{\mathcal{K}e^{-2\rho\sqrt{n+1}}}{(n+1)^{\alpha/2+1/4}},
\end{equation}
where $\mathcal{K}\cong\widehat{V}_{\alpha}/(2\sqrt{\pi})$ for $n\gg1$. 

\item[\rm(ii)] The error of $\Pi_n^{\mathrm{F}}f$ in the weighted norm satisfies
\begin{equation}\label{eq:LagErrorW}
\|f - \Pi_n^{\mathrm{F}}f\|_{\varpi_{\alpha}} \leq
\mathcal{K}e^{-2\rho\sqrt{n+1}},
\end{equation}
where $\mathcal{K}\cong\widehat{V}_{\alpha}/(2\sqrt{2\pi\rho})$ for $n\gg1$.

\item[\rm(iii)] For large $n$, the error of $\Pi_n^{\mathrm{F}}f$ in the maximum norm satisfies
\begin{align}\label{eq:LagErrorM}
\|f - \Pi_n^{\mathrm{F}}f\|_{\infty} &\leq \mathcal{K}
(n+1)^{|\alpha|/2+1/4}e^{-2\rho\sqrt{n+1}},
\end{align}
where $\|f\|_{\infty}=\max_{x\in\mathbb{R}_{+}}|f(x)|$ and
\[
\mathcal{K} \cong
\frac{\widehat{V}_{\alpha}}{2\sqrt{\pi}\rho}\left\{
\begin{array}{ll}
{\displaystyle \frac{1}{\Gamma(\alpha+1)}},   &  \alpha\geq0,  \\[12pt]
{\displaystyle 2},       &  \alpha\in(-1,0).
\end{array}
\right.
\]
\end{itemize}
\end{theorem}
\begin{proof}
Observe that the Laguerre coefficient $b_n^{(\alpha)}$ of $f(x)$ is
exactly the Laguerre coefficient $a_n^{(\alpha)}$ of $e^{x/2}f(x)$
defined in \eqref{def:LagProj}. Hence \eqref{eq:akF} follows
directly from Theorem \ref{thm:LagPolyBound}. As for
\eqref{eq:LagErrorW}, using \eqref{eq:akF} and \eqref{eq:gammaAsym},
we obtain
\begin{align}
\|f - \Pi_n^{\mathrm{F}}f\|_{\varpi_{\alpha}}^{2} &=
\sum_{k=n+1}^{\infty} (b_{k}^{(\alpha)})^2 \gamma_k^{(\alpha)} \leq
\mathcal{K}^2 \sum_{k=n+1}^{\infty}
\frac{e^{-4\rho\sqrt{k+1}}}{\sqrt{k+1}} \leq
\frac{\mathcal{K}^2}{2\rho} e^{-4\rho\sqrt{n+1}},  \nonumber
\end{align}
where $\mathcal{K}\cong\widehat{V}_{\alpha}/(2\sqrt{\pi})$. 
The desired result \eqref{eq:LagErrorW} then follows by taking a
square root on both sides of the above inequality. Finally, we
consider the proof of \eqref{eq:LagErrorM}. Here we only consider
the case of $\alpha\geq0$ since the case of $\alpha\in(-1,0)$ can be
proved in a similar way. By \eqref{eq:gammaAsym},
\eqref{eq:LagFunBound} and \eqref{eq:akF}, we have
\begin{align}
\|f - \Pi_n^{\mathrm{F}}f\|_{\infty} \leq \sum_{k=n+1}^{\infty}
|b_k^{(\alpha)} \widehat{L}_k^{(\alpha)}(x)| &\leq
\mathcal{K}\sum_{k=n+2}^{\infty} k^{\alpha/2-1/4}
e^{-2\rho\sqrt{k}}. \nonumber
\end{align}
where
$\mathcal{K}\cong\widehat{V}_{\alpha}/(2\sqrt{\pi}\Gamma(\alpha+1))$. 
Let $\Upsilon$ denote the sum in the last inequality, for large $n$,
it can be bounded by
\begin{align}
\Upsilon \leq \int_{n+1}^{\infty} x^{\alpha/2-1/4}
e^{-2\rho\sqrt{x}} \mathrm{d}x &=
\frac{\Gamma(\alpha+3/2,2\rho\sqrt{n+1})}{2^{\alpha+1/2}
\rho^{\alpha+3/2}}, \nonumber
\end{align}
where $\Gamma(a,z)$ is the incomplete gamma function
\cite[Chapter~8]{Olver2010}. Furthermore, recall that
$\Gamma(a,z)=z^{a-1}e^{-z}(1+O(z^{-1}))$ for large $z$ (see, e.g.,
\cite[Equation~(8.11.2)]{Olver2010}), the desired result
\eqref{eq:LagErrorM} follows from the last two inequalities. This
ends the proof.
\end{proof}

\begin{figure}[ht]
\includegraphics[width=0.49\textwidth,height=0.42\textwidth]{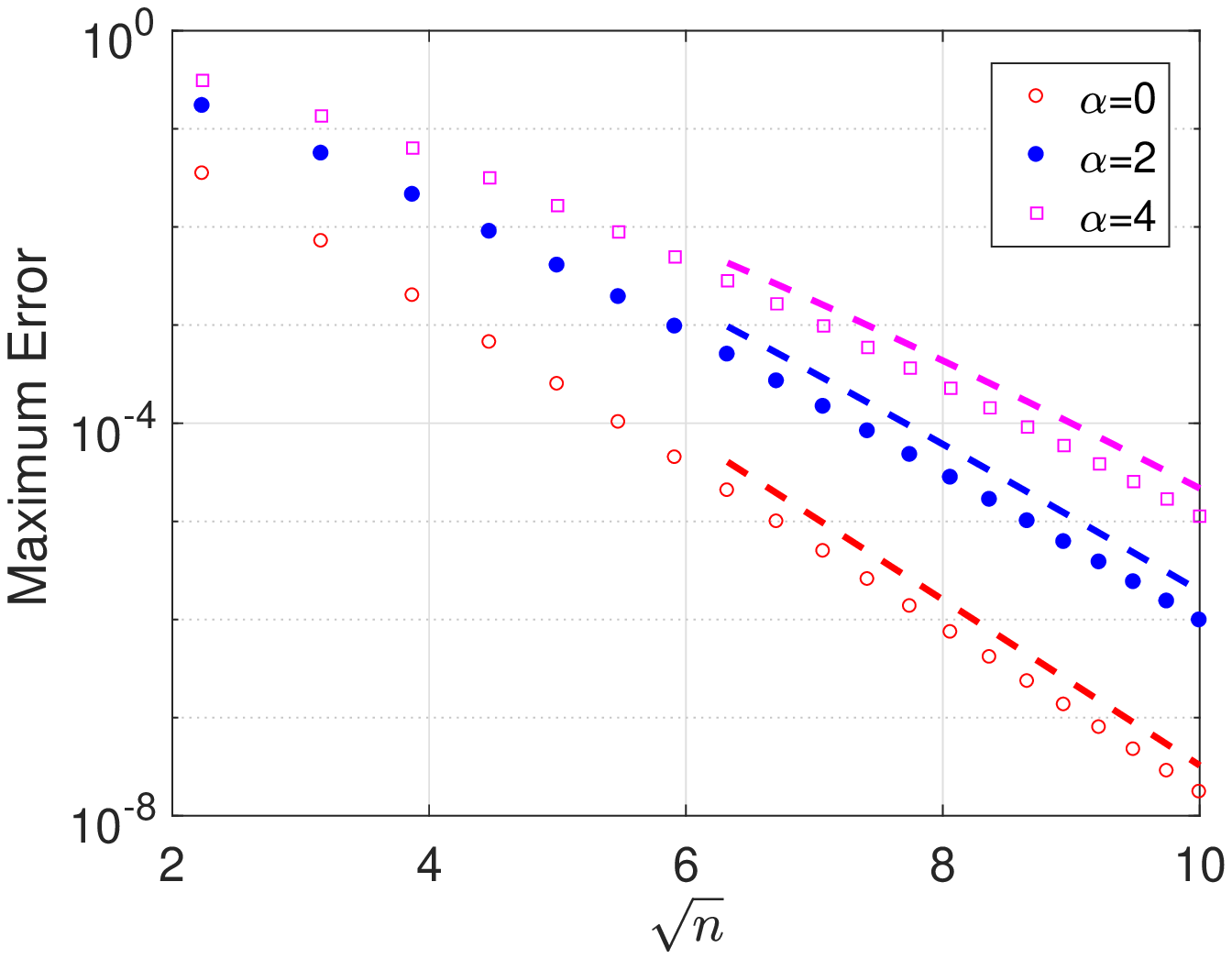}~
\includegraphics[width=0.49\textwidth,height=0.42\textwidth]{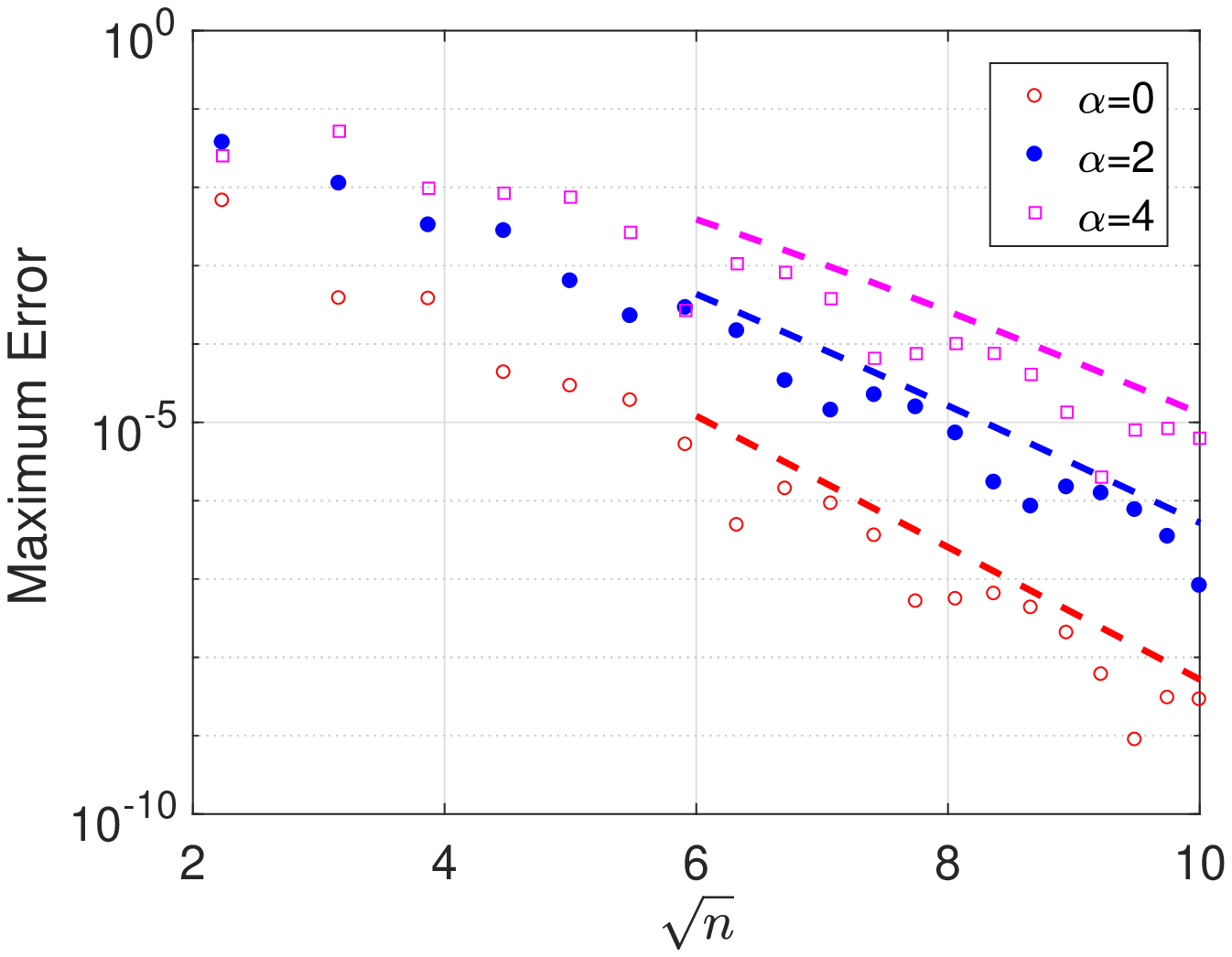}
\caption{The maximum errors of $\Pi_n^{\mathrm{F}}f$ and the
predicted convergence rates (dashed) for $f(x)=e^{-x/2}/(1+x)$
(left) and $f(x)=e^{-2x/3}/(x^2+4)$ (right).}\label{fig:LagFunApp}
\end{figure}

Figure \ref{fig:LagFunApp} shows the maximum errors of
$\Pi_n^{\mathrm{F}}f$ for $f(x)=e^{-x/2}/(1+x)$ and
$f(x)=e^{-2x/3}/(x^2+4)$. It is easily seen that the former has a
simple pole at $z=-1$ and $|f(z)|=O(|z|^{-1}e^{-\Re(z)/2})$ and the
latter has a pair of complex conjugate poles at $z=\pm2\mathrm{i}$
and $|f(z)|=O(|z|^{-2}e^{-2\Re(z)/3})$ as $z\rightarrow\infty$, and
both satisfy the assumptions of Theorem \ref{thm:LagFunBound}.
Hence, the predicted rates of convergence of $\Pi_n^{\mathrm{F}}f$
for these two functions are $O(n^{|\alpha|/2+1/4}e^{-2\sqrt{n}})$.
From Figure \ref{fig:LagFunApp} we can see that the predicted rates
of convergence are in good agreement with the actual rates of
convergence.

\section{Convergence rate analysis of Laguerre interpolation}\label{sec:LagInterp}
In this section we consider the convergence rate analysis of
Laguerre interpolation. Before embarking upon this, we state a
useful lemma about the remainder of polynomial interpolation on
$\mathbb{R}_{+}$.
\begin{lemma}\label{lem:Remainder}
Let $f$ be analytic inside and on the parabola $P_{\rho}$ with
$\rho>0$ and let $|f(z)|\leq\mathcal{K}|z|^{\beta}$ for some
$\beta<1/2$ as $|z|\rightarrow\infty$ in $D_{\rho}$. Moreover, let
$0\leq{x}_0<x_1<\cdots<x_n<\infty$ be a set of distinct points on
$\mathbb{R}_{+}$ and let $p_n$ be the unique polynomial of degree
$n$ which interpolates $f$ at these points. Then, for
$x\in\mathbb{R}_{+}$, we have
\begin{equation}\label{eq:LagInterp}
f(x) - p_n(x) = \frac{1}{2\pi\mathrm{i}} \int_{P_{\rho}}
\frac{\phi(x)f(z)}{\phi(z)(z-x)} \mathrm{d}z,
\end{equation}
where $\phi(z)=c\prod_{j=0}^{n}(z-x_j)$ and $c$ is an arbitrary
nonzero constant.
\end{lemma}
\begin{proof}
Let $\Omega$ and $\Gamma$ be defined as in the proof of Theorem
\ref{thm:LagParab} and let $\eta>x_n$ so that $\Omega$ contains all
interpolation points. For $x\in[0,\eta-\varepsilon]$ where
$\varepsilon>0$ is small, using Hermite's contour integral for the
remainder of polynomial interpolation
\cite[Theorem~3.6.1]{Davis1975}, we obtain
\begin{equation}
f(x) - p_n(x) = \frac{1}{2\pi\mathrm{i}} \int_{\Gamma+\mathcal{V}}
\frac{\phi(x)f(z)}{\phi(z)(z-x)}\mathrm{d}z,  \nonumber
\end{equation}
where $\mathcal{V}$ denotes the portion of the vertical line
$\Re(z)=\eta$ lying inside the region $D_{\rho}$. It is easily
verified that $\min_{z\in\mathcal{V}}|\phi(z)|=|\phi(\eta)|$, and
thus
\begin{align}
\left|\int_{\mathcal{V}} \frac{\phi(x)f(z)}{\phi(z)(z-x)}
\mathrm{d}z \right| \leq \left| \frac{\phi(x)}{\phi(\eta)} \right|
\int_{\mathcal{V}} \left| \frac{f(z)}{z-x} \right| \mathrm{d}s \leq
\int_{\mathcal{V}} \left| \frac{f(z)}{z-x} \right| \mathrm{d}s
\rightarrow0, \quad \eta\rightarrow\infty, \nonumber
\end{align}
where we have used the fact that $|\phi(x)|\leq|\phi(\eta)|$ for
$x\in[0,\eta-\varepsilon]$ and for large $\eta$ in the second step
and the assumption $|f(z)|\leq\mathcal{K}|z|^{\beta}$ for some
$\beta<1/2$ as $|z|\rightarrow\infty$ in $D_{\rho}$ in the last
step. Since $\Gamma\rightarrow{P}_{\rho}$ as
$\eta\rightarrow\infty$, the desired result \eqref{eq:LagInterp}
follows. This ends the proof.
\end{proof}

In the following we will consider the error estimate of Laguerre
interpolation. We divide our discussion into two cases: one is
interpolation using polynomials and the other is interpolation using
GLFs. We start by considering the former case.
\begin{theorem}\label{thm:Laginterp}
Under the assumptions of Lemma \ref{lem:Remainder} and assuming that
the interpolation points are the Laguerre points, i.e.,
$\{x_j\}_{j=0}^{n}$ are the zeros of $L_{n+1}^{(\alpha)}(x)$, the
error of the corresponding interpolation polynomial satisfies
\begin{equation}\label{eq:LaginterpI}
\|f - p_n\|_{\omega_{\alpha}} \leq \mathcal{K} n^{1/4}
e^{-2\rho\sqrt{n}},
\end{equation}
where $\mathcal{K}\cong V_{\alpha}/(\sqrt{\pi}\rho^2)$ for $n\gg1$. 
If the interpolation points are the Laguerre-Radau points, i.e.,
$\{x_j\}_{j=0}^{n}$ are the zeros of $xL_{n}^{(\alpha+1)}(x)$, the
error of the corresponding interpolation polynomial satisfies
\begin{equation}\label{eq:LaginterpII}
\|f - p_n\|_{\omega_{\alpha}} \leq \mathcal{K} n^{3/4}
e^{-2\rho\sqrt{n}},
\end{equation}
where $\mathcal{K}\cong{V}_{\alpha-1}\sqrt{2/\pi}/\rho^2$ for $n\gg1$. 
\end{theorem}
\begin{proof}
From Lemma \ref{lem:Remainder} it follows that
\begin{align}\label{eq:LagInterpBound}
\|f-p_n\|_{\omega_{\alpha}} &\leq \frac{1}{2\pi}
\sqrt{\int_{\mathbb{R}_{+}} \omega_{\alpha}(x) \left[
\int_{P_{\rho}} \frac{|\phi(x)f(z)|}{|\phi(z)(z-x)|} \mathrm{d}s
\right]^2 \mathrm{d}x } \leq
\frac{\|\phi\|_{\omega_{\alpha}}}{2\pi\rho^2} \int_{P_{\rho}}
\frac{|f(z)|}{|\phi(z)|} \mathrm{d}s,
\end{align}
where we have used the fact that $\min|z-x|=\rho^2$ for
$z\in{P}_{\rho}$ and $x\in\mathbb{R}_{+}$ in the last step. We first
consider the case where the interpolation points are the Laguerre
points, i.e., $\phi(x)=L_{n+1}^{(\alpha)}(x)$, then using
\eqref{eq:gammaAsym} and Perron's formula for GLPs
\cite[Theorem~8.22.3]{Szego1975}, we have
\begin{align}
\|\phi\|_{\omega_{\alpha}} = \gamma_{n+1}^{(\alpha)}, \quad
\int_{P_{\rho}} \frac{|f(z)|}{|\phi(z)|} \mathrm{d}s &=
\frac{e^{-2\rho\sqrt{n}}}{n^{\alpha/2-1/4}}
\left(2\sqrt{\pi}V_{\alpha}+O(n^{-1/2})\right). \nonumber
\end{align}
Combining these two results with \eqref{eq:LagInterpBound} gives
\eqref{eq:LaginterpI}. Next, we consider the case where the
interpolation points are the Laguerre-Radau points, i.e.,
$\phi(x)=xL_{n}^{(\alpha+1)}(x)$. By
\cite[Equation~(18.9.14)]{Olver2010}, we have that
\begin{align}
\|\phi\|_{\omega_{\alpha}} = \sqrt{(n+\alpha+1)^2
\gamma_n^{(\alpha)} + (n+1)^2\gamma_{n+1}^{(\alpha)}}. \nonumber
\end{align}
Moreover, using Perron's formula for GLPs
\cite[Theorem~8.22.3]{Szego1975} again, we get
\begin{align}
\int_{P_{\rho}} \frac{|f(z)|}{|\phi(z)|} \mathrm{d}s =
\frac{e^{-2\rho\sqrt{n}}}{n^{\alpha/2+1/4}} \left(
2\sqrt{\pi}V_{\alpha-1} + O(n^{-1/2})\right). \nonumber
\end{align}
Combining the above two results with \eqref{eq:LagInterpBound} gives
\eqref{eq:LaginterpII}. This ends the proof.
\end{proof}

\begin{figure}[ht]
\centering
\includegraphics[width=0.48\textwidth,height=0.42\textwidth]{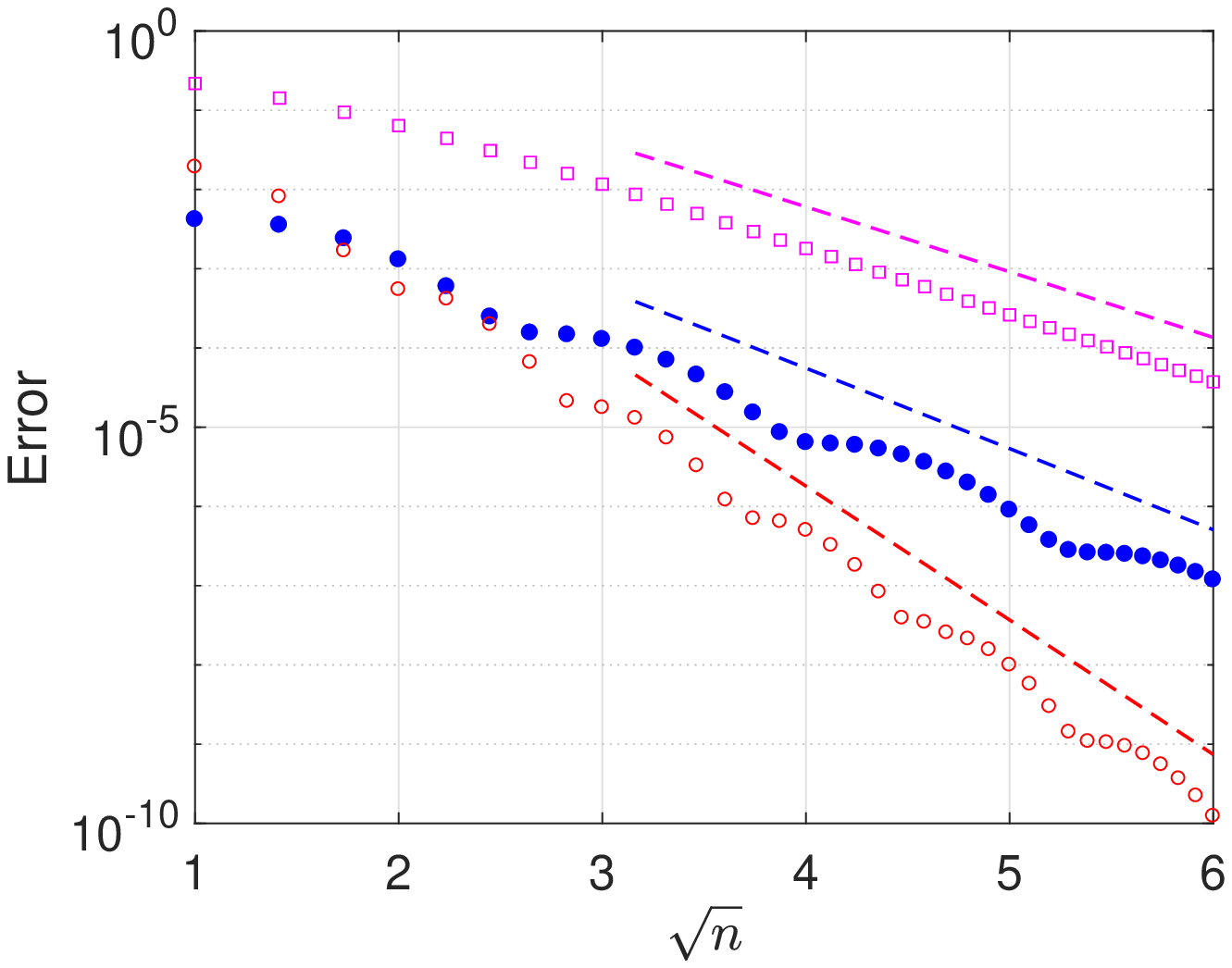}\quad
\includegraphics[width=0.48\textwidth,height=0.42\textwidth]{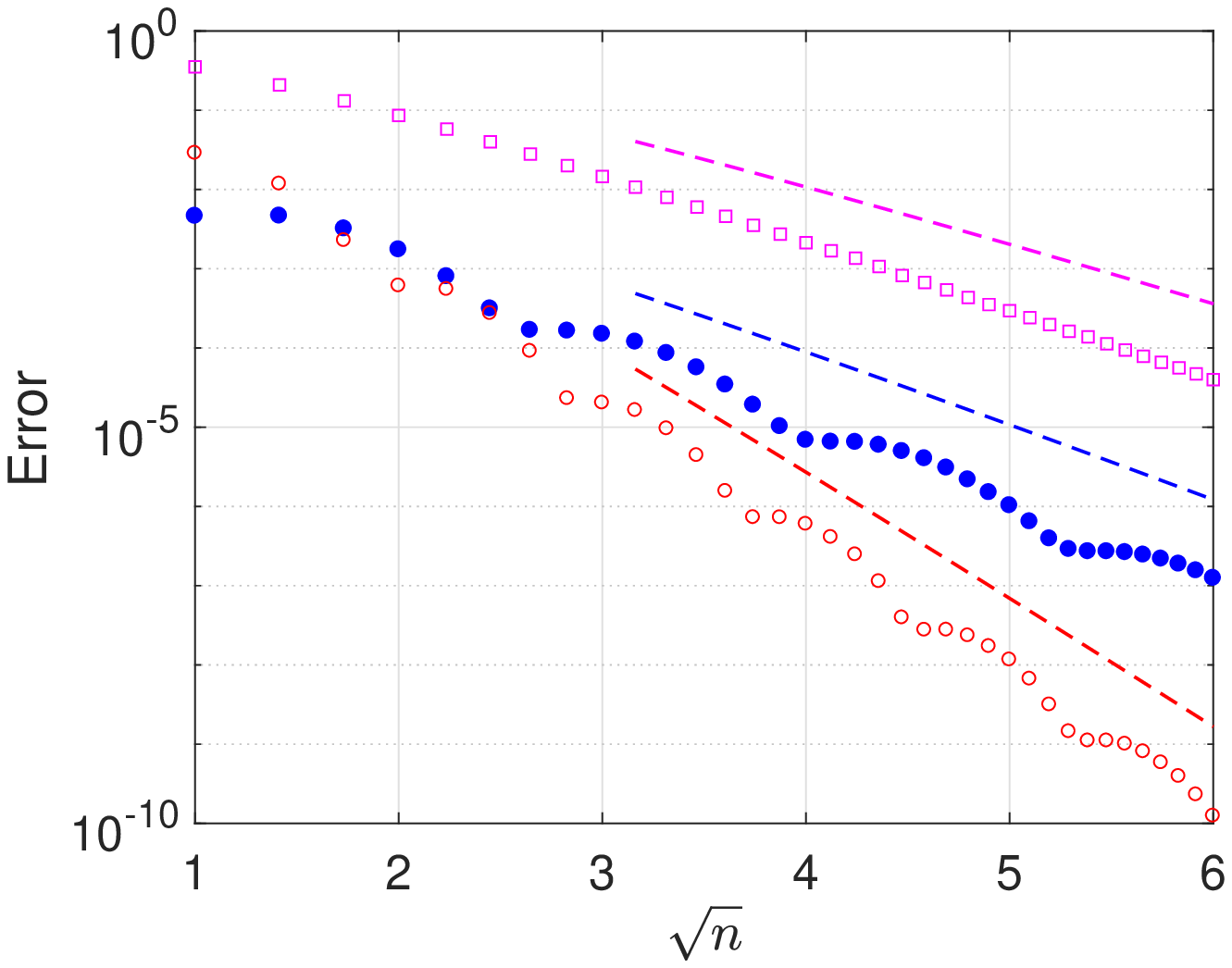}
\caption{The errors of the Laguerre interpolants (left) and
Laguerre-Radau interpolants (right) for $f(x)=e^{-x}/(x+1)$ (boxes),
$f(x)=1/(x^2+9)$ (dots) and $f(x)=\mathrm{sech}(\pi{x}/16)$
(circles). The dashed lines denote the predicted convergence rates.
}\label{fig:LagInterp}
\end{figure}

In Figure \ref{fig:LagInterp} we display the errors of Laguerre and
Laguerre-Radau interpolants for three test functions. In our
computations, the Laguerre and Laguerre-Radau interpolants are
evaluated by using their barycentric form
\begin{equation}
p_n(x) = \frac{\displaystyle \sum_{j=0}^{n}\frac{\lambda_j}{x-x_j}
f(x_j)}{\displaystyle \sum_{j=0}^{n}\frac{\lambda_j}{x-x_j}},
\end{equation}
where $\{\lambda_j\}_{j=0}^{n}$ are the barycentric weights
associated with the points $\{x_j\}_{j=0}^{n}$. As shown in
\cite[Corollary~2.4]{Wang2014} and \cite[Theorem~3.8]{Wang2014}, the
barycentric weights of Laguerre and Laguerre-Radau points can be
calculated directly from the nodes and weights of Gauss-Laguerre and
Gauss-Laguerre-Radau quadrature rules, respectively, and fast
algorithms with a cost of $O(n)$ for these classical Gaussian
quadrature rules are available (see, e.g., \cite{Opsomer2022}). We
can see from Figure \ref{fig:LagInterp} that the errors of Laguerre
interpolants are consistent with our analysis. As for Laguerre-Radau
interpolants, numerical results show that their convergence rate is
still $O(n^{1/4}e^{-2\rho\sqrt{n}})$, which is slightly faster than
the predicted one.

Finally, we consider the error estimates of Laguerre interpolation
using GLFs. Specifically, let $q_n\in\mathbb{Q}_n$ such that
\begin{equation}\label{eq:ModLagInterp}
q_n(x_j) = f(x_j), \quad j=0,\ldots,n,
\end{equation}
where $\{x_j\}_{j=0}^{n}$ are the Laguerre or the Laguerre-Radau
points. Error estimates for $q_n$ in the maximum norm are given
below.
\begin{theorem}
Let $f$ be analytic inside and on the parabola $P_{\rho}$ with
$\rho>0$ and let $|f(z)|\leq\mathcal{K}|z|^{\beta}e^{-\Re(z)/2}$ for
some $\beta<1/2$ as $|z|\rightarrow\infty$ in $D_{\rho}$. If the
points $\{x_j\}_{j=0}^{n}$ are the Laguerre points, then
\begin{equation}\label{eq:LagfunInterp1}
\|f - q_n\|_{\infty} \leq \mathcal{K}
n^{|\alpha|/2+1/4}e^{-2\rho\sqrt{n}},
\end{equation}
where
$\mathcal{K}\cong\widehat{V}_{\alpha}/(\sqrt{\pi}\rho^2\Gamma(\alpha+1))$
for $\alpha\geq0$ and
$\mathcal{K}\cong2\widehat{V}_{\alpha}/(\sqrt{\pi}\rho^2)$ for
$\alpha\in(-1,0)$. If the points $\{x_j\}_{j=0}^{n}$ are the
Laguerre-Radau points, then
\begin{equation}\label{eq:LagfunInterp2}
\|f - q_n\|_{\infty} \leq \mathcal{K}
n^{|\alpha|/2+3/4}e^{-2\rho\sqrt{n}},
\end{equation}
where
$\mathcal{K}\cong2\widehat{V}_{\alpha-1}/(\sqrt{\pi}\rho^2\Gamma(\alpha+1))$
for $\alpha\geq0$ and
$\mathcal{K}\cong2\widehat{V}_{\alpha-1}/(\sqrt{\pi}\rho^2)$ for
$\alpha\in(-1,0)$.
\end{theorem}
\begin{proof}
We only consider the proof of the case of Laguerre points. From
\eqref{eq:ModLagInterp} it is easily verified that $e^{x/2}q_n(x)$
is the unique polynomial of degree $n$ which interpolates
$e^{x/2}f(x)$ at the Laguerre points. Hence, by Lemma
\ref{lem:Remainder} we have that
\begin{equation}
e^{x/2}f(x) - e^{x/2}q_n(x) = \frac{1}{2\pi\mathrm{i}}
\int_{P_{\rho}} \frac{L_{n+1}^{(\alpha)}(x)
e^{z/2}f(z)}{L_{n+1}^{(\alpha)}(z) (z-x)} \mathrm{d}z, \nonumber
\end{equation}
or equivalently,
\begin{equation}
f(x) - q_n(x) = \frac{1}{2\pi\mathrm{i}} \int_{P_{\rho}}
\frac{\widehat{L}_{n+1}^{(\alpha)}(x)
e^{z/2}f(z)}{L_{n+1}^{(\alpha)}(z)(z-x)} \mathrm{d}z. \nonumber
\end{equation}
The desired result \eqref{eq:LagfunInterp1} follows by combining
\eqref{eq:LagFunBound} and the Perron's formula for GLPs
\cite[Theorem~8.22.3]{Szego1975}. This ends the proof.
\end{proof}

\section{Laguerre approximations for entire functions}\label{sec:Fast}
\begin{figure}[ht]
\centering
\includegraphics[width=0.48\textwidth,height=0.42\textwidth]{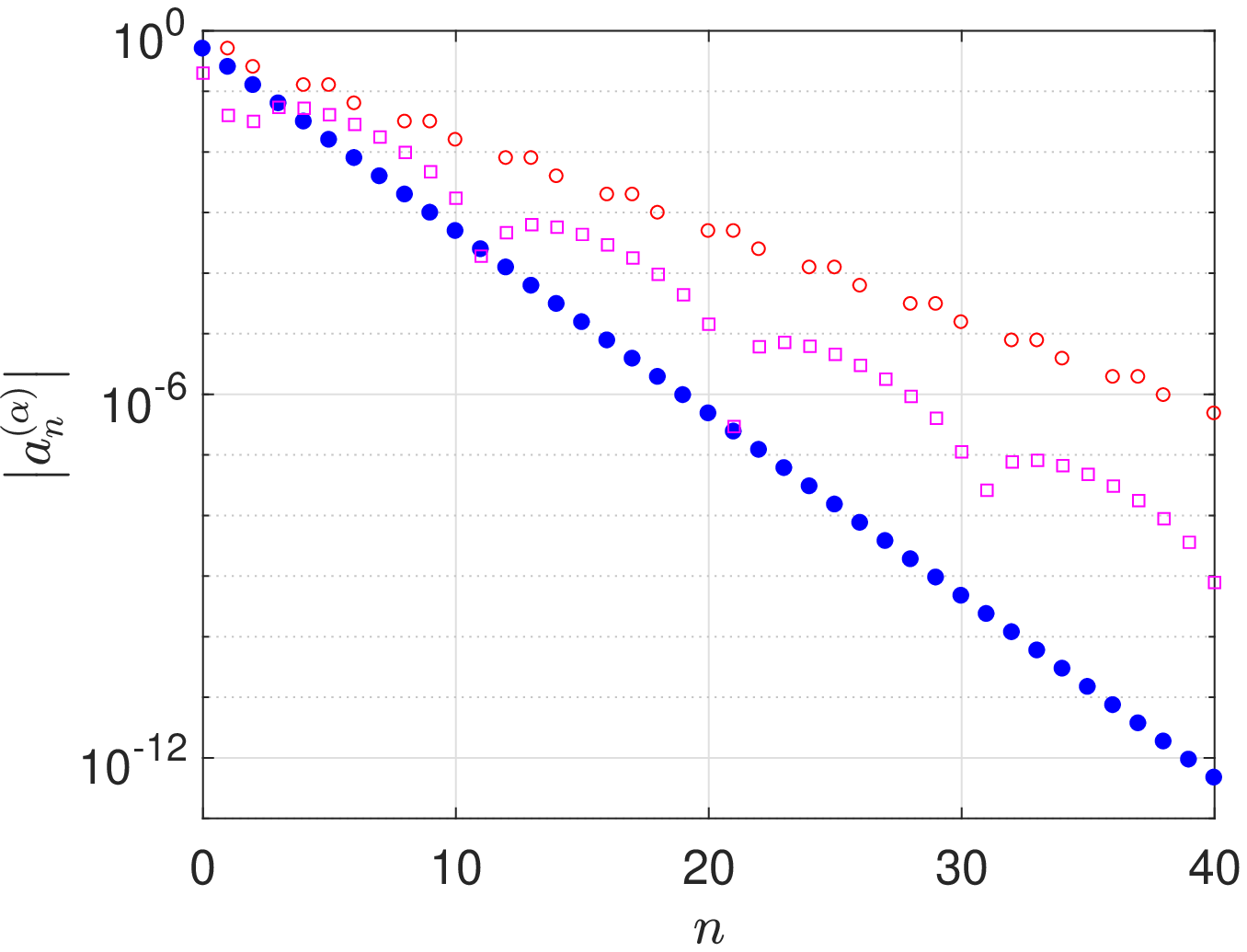}\quad
\includegraphics[width=0.48\textwidth,height=0.42\textwidth]{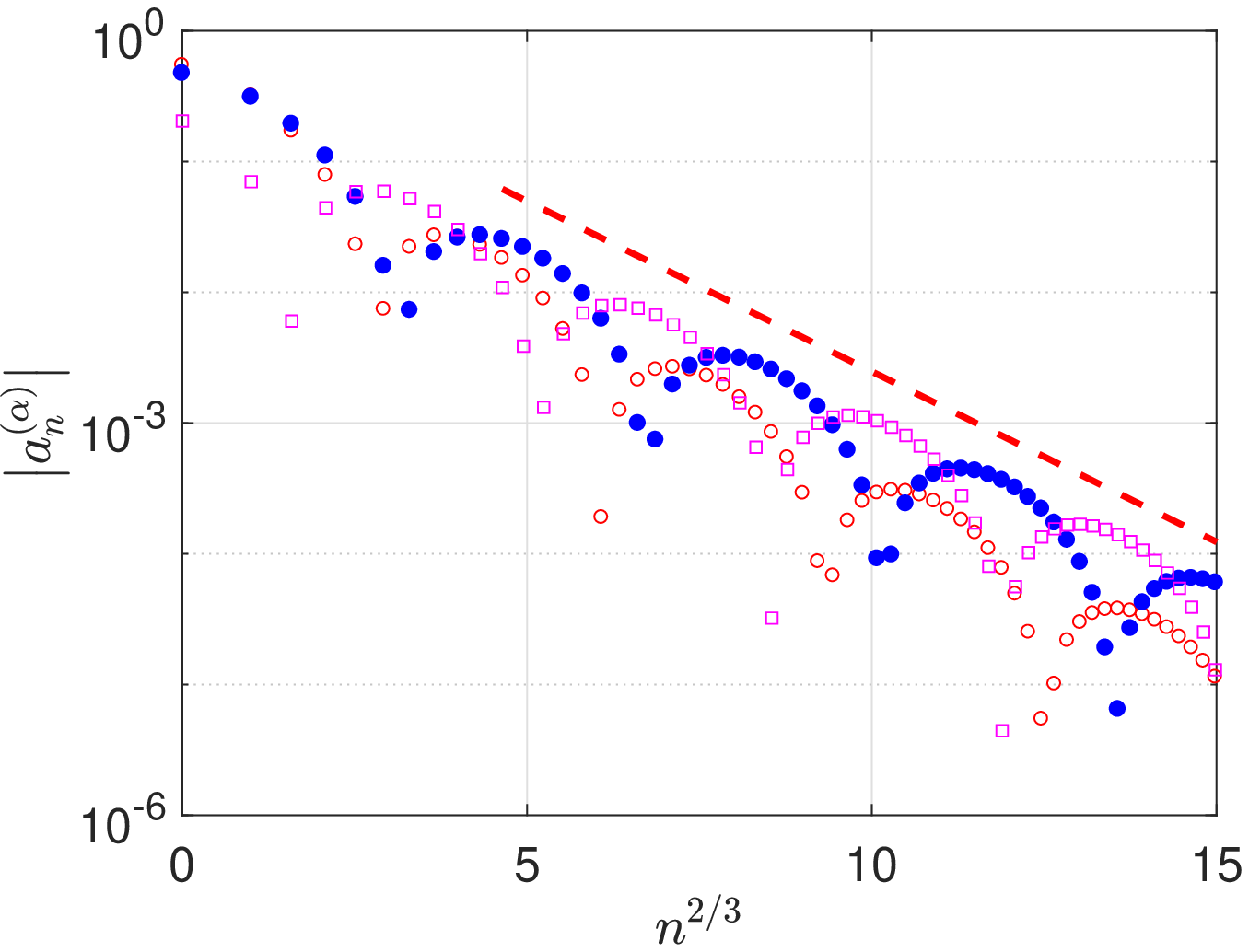}
\caption{The magnitude of the Laguerre coefficients of $f(x)$. The
left panel shows $f(x)=\cos(x)$ (circles), $f(x)=e^{-x}\sin(x)$
(boxes) and $f(x)=e^{-x}$ (dots) and the right panel shows
$f(x)=e^{-x^2}$ (circles), $f(x)=e^{-x^2}\sin(x)$ (boxes) and
$f(x)=e^{-x^2}\cos(x)$ (dots). The dashed line in the right panel is
$e^{-\kappa{n}^{2/3}}$ with $\kappa=3/5$.}\label{fig:LagExp}
\end{figure}

In the above two sections we have proved the root-exponential
convergence of Laguerre approximations for analytic functions. An
intriguing question is whether and when it is possible to improve
this convergence. Indeed, the root-exponential convergence might be
improved in the context of entire functions. To show this, consider
a set of test functions $f(x)=\cos(x),e^{-x}\sin(x),e^{-x}$ which
are all entire functions of exponential type. The left panel of
Figure \ref{fig:LagExp} shows the magnitudes of the Laguerre
coefficients of these functions, from which it is easily seen that
the Laguerre coefficients decay at exponential rates. 
We also consider another set of test functions
$f(x)=e^{-x^2},e^{-x^2}\sin(x),e^{-x^2}\cos(x)$ which are all entire
functions of Gaussian type. The right panel of Figure
\ref{fig:LagExp} shows the magnitudes of the Laguerre coefficients
of these functions, from which we see that the Laguerre coefficients
behave like $a_n^{(\alpha)}=O(e^{-\kappa n^{2/3}})$ for some
$\kappa>0$. Note that the rate of convergence of Laguerre
approximations is the same as that of decay of the Laguerre
coefficients.

We conclude from the above observations that the root-exponential
convergence of Laguerre approximations might be improved when the
underlying functions are entire. Furthermore, the precise rate of
convergence of Laguerre approximations for entire functions depends
on the behavior of the underlying function at infinity and their
analysis is a delicate issue. We omit further discussion here.

\section{Applications}\label{sec:Extension}
In this section we extend our study to several important applications,
including Laguerre spectral differentiations, Gauss-Laguerre
quadrature rules, the scaling factors and the Weeks method for the numerical inversion of
the Laplace transform.

\subsection{Laguerre spectral differentiations}
We consider error estimates of Laguerre spectral differentiations.
The main results are stated in the following theorem.
\begin{theorem}
Let $\Pi_n^{\mathrm{P}}f$ and $\Pi_n^{\mathrm{F}}f$ be defined in
\eqref{def:LagProj} and \eqref{def:LagProj2}, respectively. For
$m\in\mathbb{N}$, the following statements hold.
\begin{itemize}
\item[\rm(i)] Under the assumptions of Theorem \ref{thm:LagPolyBound},
then for each $n\geq m$,
\begin{align}\label{eq:LSDI}
\|f^{(m)}-(\Pi_n^{\mathrm{P}}f)^{(m)}\|_{\omega_{\alpha+m}} \leq
\mathcal{K} n^{m/2} e^{-2\rho\sqrt{n}},
\end{align}
where $\mathcal{K}\cong V_{\alpha}/(2\sqrt{2\pi\rho})$ for $n\gg1$. 

\item[\rm(ii)] Under the assumptions of Theorem \ref{thm:LagFunBound},
then for each $n\geq m$,
\begin{align}\label{eq:LSDII}
\|f^{(m)} - (\Pi_n^{\mathrm{F}}f)^{(m)}\|_{\infty} \leq
\mathcal{K}n^{m+\alpha/2+1/4} e^{-2\rho\sqrt{n}}.
\end{align}
where
$\mathcal{K}\cong\widehat{V}_{\alpha}/(2\sqrt{\pi}\rho\Gamma(\alpha+m+1))$ for $n\gg1$. 
\end{itemize}
\end{theorem}
\begin{proof}
From \cite[Equation~(18.9.23)]{Olver2010} it follows that
\begin{align}
\|f^{(m)}-(\Pi_n^{\mathrm{P}}f)^{(m)}\|_{\omega_{\alpha+m}}^2 &=
\sum_{k=n+1}^{\infty} (a_k^{(\alpha)})^2 \gamma_{k-m}^{(\alpha+m)}.
\nonumber
\end{align}
Combining this with Theorem \ref{thm:LagPolyBound} gives
\eqref{eq:LSDI}. As for \eqref{eq:LSDII}, using
\cite[Equation~(18.9.23)]{Olver2010} again and \eqref{eq:gammaAsym}
and \eqref{eq:LagFunBound}, one can show that
$\|\partial_{x}^{m}\widehat{L}_n^{(\alpha)}\|_{\infty}\leq
\mathcal{K}n^{\alpha+m}$, where
$\mathcal{K}\cong1/\Gamma(\alpha+m+1)$. Note that
\begin{align}
\|f^{(m)} - (\Pi_n^{\mathrm{F}}f)^{(m)}\|_{\infty} \leq
\sum_{k=n+1}^{\infty} |b_k^{(\alpha)}|
\|\partial_{x}^{m}\widehat{L}_n^{(\alpha)}\|_{\infty}. \nonumber
\end{align}
The desired result \eqref{eq:LSDII} follows from Theorem
\ref{thm:LagFunBound}. This ends the proof.
\end{proof}

\begin{remark}
For Laguerre spectral differentiations based on interpolation
methods, their root-exponential convergence can also be analyzed by
using contour integral representations presented in section
\ref{sec:LagInterp}.
\end{remark}

\subsection{Gauss-Laguerre quadrature rules}
Consider integrals of the form
\begin{equation}
I(f) = \int_{\mathbb{R}_{+}}\omega_{\alpha}(x)f(x)\mathrm{d}x,
\end{equation}
and their interpolatory quadrature rules
\begin{equation}
Q_n(f) = \sum_{k=0}^{n} w_k f(x_k).
\end{equation}
It is well-known that, among all interpolatory quadrature rules,
Gaussian quadrature rules are optimal in the sense that they have
the highest degree of exactness when using the same number of
quadrature nodes. More specifically, an $(n+1)$-point Gauss-Laguerre
quadrature rule has degree of exactness $2n+1$ and an $(n+1)$-point
Gauss-Laguerre-Radau quadrature with a preassigned node $x_0=0$ has
degree of exactness $2n$. Below, for notational simplicity, we
denote by $Q_n^{G}(f)$ the $(n+1)$-point Gauss-Laguerre quadrature
and $Q_n^{R}(f)$ the $(n+1)$-point Gauss-Laguerre-Radau quadrature.
Our main result is stated in the following theorem.
\begin{theorem}\label{thm:GLQuad}
Under the same assumptions of Theorem \ref{thm:LagParab}, it holds
that
\begin{align}\label{eq:GLQuad}
\left.
\begin{array}{ll}
|I(f) - Q_n^{G}(f)| \\[8pt]
|I(f) - Q_n^{R}(f)|
\end{array}
\right\} \leq \mathcal{K}e^{-4\rho\sqrt{n}},
\end{align}
where
$\mathcal{K}\cong\int_{P_{\rho}}|e^{-z}(-z)^{\alpha+1/2}f(z)|\mathrm{d}s$ for $n\gg1$.
\end{theorem}
\begin{proof}
We start by considering the case of Gauss-Laguerre quadrature.
Recall that Gauss-Laguerre is an interpolatory quadrature rule, by
Lemma \ref{lem:Remainder} we have
\begin{align}
I(f) - Q_n^{G}(f) &= \int_{\mathbb{R}_{+}}
\frac{\omega_{\alpha}(x)}{2\pi\mathrm{i}} \int_{P_{\rho}}
\frac{L_{n+1}^{(\alpha)}(x)f(z)}{L_{n+1}^{(\alpha)}(z)(z-x)}
\mathrm{d}z \mathrm{d}x = \int_{P_{\rho}}
\frac{\Phi_{n+1}^{(\alpha)}(z)}{{L_{n+1}^{(\alpha)}(z)}} f(z)
\mathrm{d}z, \nonumber
\end{align}
where we have interchanged the order of integration in the last
step. The desired result \eqref{eq:GLQuad} then follows by applying
\eqref{eq:PhiAsym} and the Perron's formula
\cite[Theorem~8.22.3]{Szego1975} to the last equation. As for the
case of Gauss-Laguerre-Radau quadrature, it follows from
\cite[Equation~(18.9.14)]{Olver2010} that
\begin{align}
I(f) - Q_n^{R}(f) &= \int_{\mathbb{R}_{+}}
\frac{\omega_{\alpha}(x)}{2\pi\mathrm{i}} \int_{P_{\rho}}
\frac{((n+\alpha+1)L_n^{(\alpha)}(x)-(n+1)L_{n+1}^{(\alpha)}(x))f(z)}{((n+\alpha+1)L_n^{(\alpha)}(z)-(n+1)L_{n+1}^{(\alpha)}(z))(z-x)}
\mathrm{d}z \mathrm{d}x \nonumber \\
&= \int_{P_{\rho}}
\frac{(n+\alpha+1)\Phi_n^{(\alpha)}(z)-(n+1)\Phi_{n+1}^{(\alpha)}(z)}{(n+\alpha+1)L_n^{(\alpha)}(z)-(n+1)L_{n+1}^{(\alpha)}(z)}
f(z) \mathrm{d}z. \nonumber
\end{align}
Applying \eqref{eq:PhiAsym} and the Perron's formula
\cite[Theorem~8.22.3]{Szego1975} to the last equation and after some
calculations gives the desired result. This ends the proof.
\end{proof}

\begin{remark}
Barrett in \cite{Barrett1961} studied the convergence properties of Gaussian quadrature rules. Under the assumptions that $f$ is analytic within and on the parabola $P_{\rho}$, except for a pair of complex conjugate simple poles, and satisfies suitable decay at infinity, Barrett derived the same root-exponential rate for the Gauss-Laguerre quadrature rule $Q_n^G(f)$. However, neither explicit proof nor the decay condition for which this rate holds were given. Moreover, our results in Theorem \ref{thm:GLQuad} are more general than Barrett's results in the sense that they hold whenever $f$ has singularities (e.g., branch points or poles) in $\mathbb{C}\setminus\mathbb{R}_{+}$, but not restricted to the case of complex conjugate simple poles.
\end{remark}

\begin{remark}
We make an interesting observation on the convergence rates of
spectral approximations and the corresponding Gaussian quadrature.
In the Laguerre case, we see that the rates of convergence of
Laguerre approximations and Gauss-Laguerre quadrature rules are
$O(e^{-2\rho\sqrt{n}})$ and $O(e^{-4\rho\sqrt{n}})$, respectively,
and thus the latter are four times faster than the former. In the
Legendre case, it has been known that the convergence rates of
Legendre approximations and Gauss-Legendre quadrature are
$O(\rho^{-n})$ and $O(\rho^{-2n})$ respectively if the underlying
function is analytic in an ellipse with foci $\pm1$ and major and
minor semiaxis lengths summing to $\rho>1$ (see
\cite[Theorem~2]{Wang2021} and \cite[Theorem~19.3]{Trefeth2013}).
Hence, the rate of convergence of Gauss-Legendre quadrature is only
twice faster than that of Legendre approximations.
\end{remark}

\begin{figure}[ht]
\centering
\includegraphics[width=0.48\textwidth,height=0.42\textwidth]{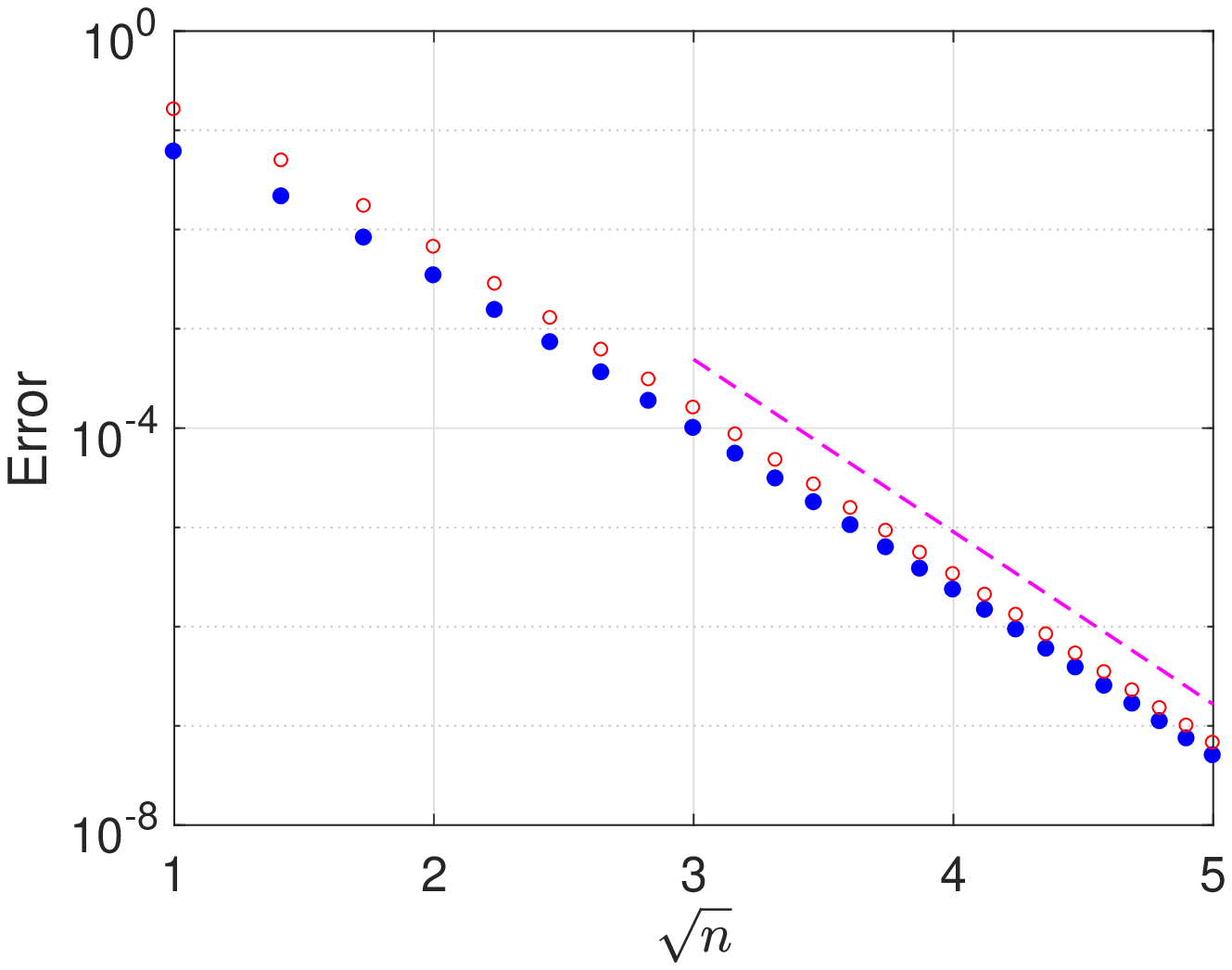}\quad
\includegraphics[width=0.48\textwidth,height=0.42\textwidth]{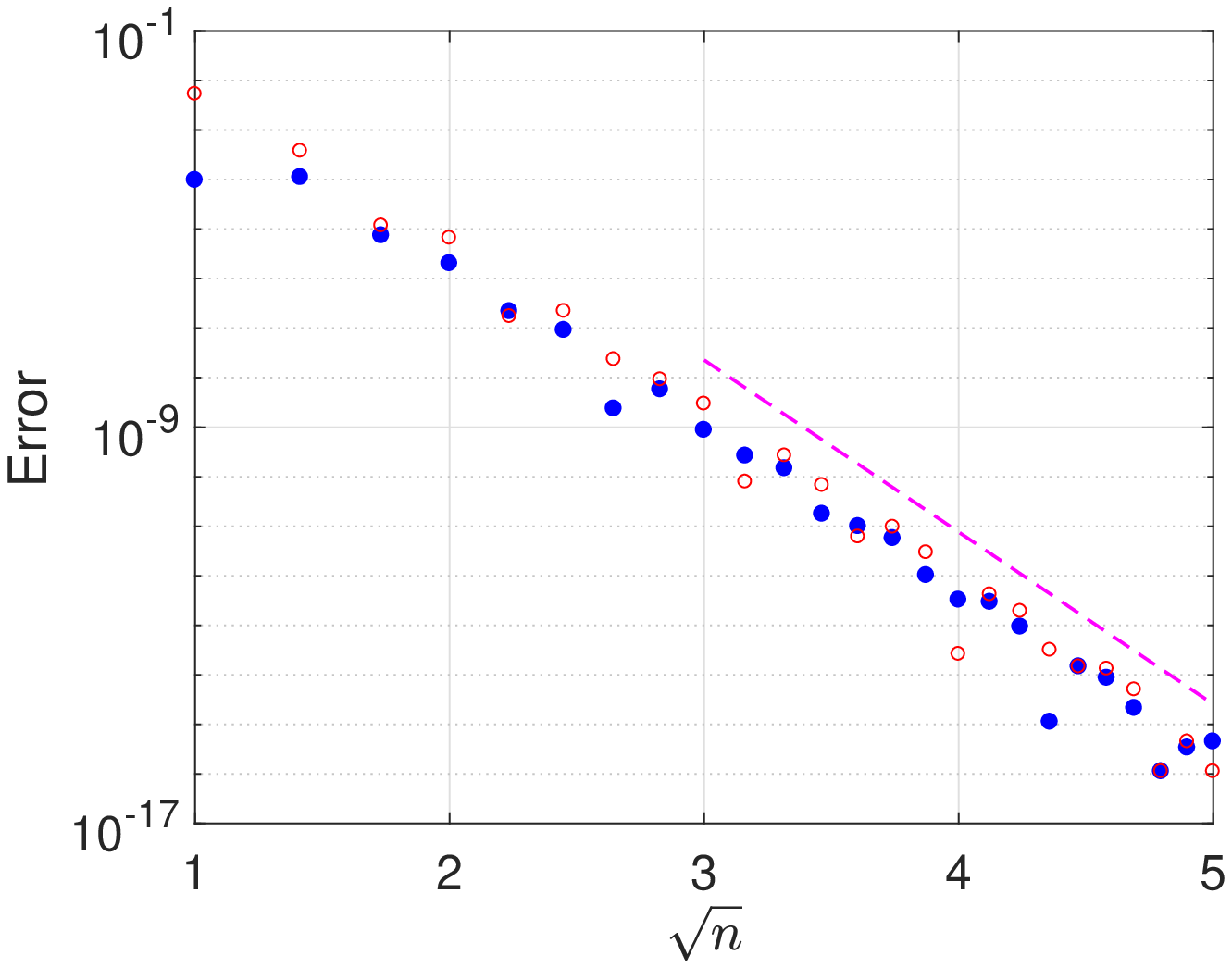}
\caption{The errors of $Q_n^{G}(f)$ (dots) and $Q_n^{R}(f)$
(circles) for $f(x)=e^{-x}/(x+1)$ and $\alpha=0$ (left) and
$f(x)=\mathrm{sech}(\pi{x}/16)$ and $\alpha=1/2$ (right). The dash
lines show the predicted convergence rates.}\label{fig:LagQuad}
\end{figure}

Figure \ref{fig:LagQuad} illustrates the errors of Gauss-Laguerre
and Gauss-Laguerre-Radau quadrature rules for $f(x)=e^{-x}/(x+1)$
and $f(x)=\mathrm{sech}(\pi{x}/16)$, which correspond to $\rho=1$
and $\rho=2$, respectively. We see that the errors of both
quadrature rules decay at the rate $O(e^{-4\rho\sqrt{n}})$ as
$n\rightarrow\infty$, which are consistent with our theoretical
analysis.

\subsection{The Weeks method for the inversion of Laplace transform}
The Laplace transform and its inverse are defined, respectively, by
\begin{align}
(\mathcal{L}f)(s) = \int_{0}^{\infty} e^{-st} f(t)\mathrm{d}t, \quad
f(t) = \frac{1}{2\pi\mathrm{i}}
\int_{\sigma-\mathrm{i}\infty}^{\sigma+\mathrm{i}\infty} e^{st}
(\mathcal{L}f)(s) \mathrm{d}s, \nonumber
\end{align}
where the straight line $\Re(s)=\sigma>\sigma_0$ is known as the
Bromwich contour and $\sigma_0$ is the real part of the right-most
singularity of $(\mathcal{L}f)(s)$, and they are widely used in
finding the solution of linear differential equations that arise in
problems of mathematical physics and engineering science. In
applications, analytical expressions for the inversion of Laplace
transforms are in general not available and numerical methods are
required to evaluate them accurately. In the past few decades, the
numerical inversion of the Laplace transform has attracted
considerable attention and we refer to \cite{Cohen2007} for a
comprehensive survey. Here we restrict ourselves to the Weeks method
due to its relation to Laguerre approximations.

The Weeks method can be formulated as follows (see, e.g.,
\cite{Giunta1988,Weeks1966,Weideman1999,Weideman2023}). First,
expanding $f$ into a Laguerre series, we have
\begin{equation}\label{eq:LagS}
f(t) = e^{\sigma{t}}\sum_{k=0}^{\infty} c_k
e^{-\nu{t}/2}L_k(\nu{t}),
\end{equation}
where $\sigma>\sigma_0$ and $\nu>0$ are scaling parameters. Taking
the Laplace transform on both sides of \eqref{eq:LagS} and making
the change of variable $s=\sigma-\nu/2+\nu/(1-z)$, we obtain
\begin{equation}\label{eq:TaylorExp}
\frac{\nu}{1-z}(\mathcal{L}f)\left(\sigma-\frac{\nu}{2}+\frac{\nu}{1-z}\right)
= \sum_{k=0}^{\infty} c_k z^k.
\end{equation}
Furthermore, making the change of variable $z=e^{\mathrm{i}\theta}$,
we obtain after some elementary calculation
\begin{equation}\label{eq:FourierExp}
\frac{\nu}{2}\left(1 +
\mathrm{i}\cot\left(\frac{\theta}{2}\right)\right)
(\mathcal{L}f)\left(\sigma +
\mathrm{i}\frac{\nu}{2}\cot\left(\frac{\theta}{2}\right)\right) =
\sum_{k=0}^{\infty} c_k e^{\mathrm{i}k\theta}.
\end{equation}
It is easily seen from \eqref{eq:TaylorExp} and
\eqref{eq:FourierExp} that the Laguerre coefficients
$\{c_k\}_{k=0}^{\infty}$ of $f$ are also the Taylor or the Fourier
coefficients of some functions that depend on $\mathcal{L}f$, and
both can serve as the starting point of the Weeks method. More
specifically, let $\{\tau_k\}_{k=0}^{\infty}$ be a sequence defined
by $\tau_0=1/2$ and $\tau_k=1$ for $k\geq1$, and let $\phi(\theta)$
denote the real part of the function on the left-hand side of
\eqref{eq:FourierExp}, Weeks in \cite{Weeks1966} suggested computing
the first $n+1$ Laguerre coefficients $\{c_k\}_{k=0}^{n}$ by
\begin{equation}\label{eq:WeekCoeff}
c_k^{w} =
\frac{2\tau_k}{n+1}\sum_{j=0}^{n}\phi(\theta_j)\cos(k\theta_{j}),
\quad \theta_j=\frac{(j+1/2)\pi}{n+1},
\end{equation}
and then approximating the inversion of the Laplace transform by
\begin{equation}\label{eq:WeekMethod}
f_n(t) = e^{\sigma{t}}\sum_{k=0}^{n} c_k^w e^{-\nu{t}/2}L_k(\nu{t}).
\end{equation}
Below we derive an error estimate of the Weeks method.
\begin{theorem}\label{thm:Weeks}
Assume that $f$ is analytic inside and on a parabola $P_{\rho}$ for
some $\rho>0$ and
$|f(z)|\leq\mathcal{K}|z|^{\beta}\exp((\sigma-\nu/2)\Re(z))$ for
some $\beta\in\mathbb{R}$ as $z\rightarrow\infty$ in the region
$D_{\rho}$. Then, for fixed $t\geq0$,
\begin{equation}\label{eq:WeeksRate}
|f(t) - f_n(t)| = O(n^{1/4}e^{-2\rho\sqrt{\nu{n}}}), \quad
n\rightarrow\infty.
\end{equation}
\end{theorem}
\begin{proof}
From \eqref{eq:WeekCoeff} we see that $\{c_k^w\}_{k=0}^{n}$ are
actually the Chebyshev coefficients of the unique polynomial of
degree $n$ which interpolates $\phi(\arccos(x))$ at the Chebyshev
points $\{\cos\theta_j\}_{j=0}^{n}$, hence
\begin{align}
|f(t)- f_n(t)| &\leq e^{\sigma{t}} \left(\sum_{k=0}^{\infty} c_k
e^{-\nu{t}/2}L_k(\nu{t}) - \sum_{k=0}^{n}
c_k^w e^{-\nu{t}/2}L_k(\nu{t}) \right) \nonumber \\
&\leq e^{\sigma{t}} \left(\sum_{k=0}^{n}
|c_k-c_k^w| + \sum_{k=n+1}^{\infty} |c_k| \right) \nonumber \\
&\leq 2e^{\sigma{t}} \sum_{k=n+1}^{\infty} |c_k|, \nonumber
\end{align}
where we have used the fact that $|e^{-x/2}L_k(x)|\leq1$ for
$x\in\mathbb{R}_{+}$ in the second step and the aliasing formula for
the Chebyshev coefficients \cite[p.~96]{Boyd2000}, i.e.,
\begin{align}
c_k^w &= c_k + \tau_k \sum_{\ell=1}^{\infty} (-1)^{\ell} \left(
c_{2\ell(n+1)-k} + c_{2\ell(n+1)+k} \right), \notag
\end{align}
in the last step.
On the other hand, from \eqref{eq:LagS} and by a simple change of
variable $x=\nu{t}$ we obtain that 
\begin{equation}
c_k = \nu\int_{0}^{\infty}
e^{-\sigma{t}}f(t)e^{-\nu{t}/2}L_k(\nu{t}) \mathrm{d}t =
\int_{0}^{\infty} g(x) e^{-x/2}L_k(x) \mathrm{d}x, \nonumber
\end{equation}
where $g(x)=e^{-\sigma{x}/\nu}f(x/\nu)$. Furthermore, under the
assumptions on $f$, it follows from Theorem \ref{thm:LagFunBound}
that $|c_k|=O(k^{-1/4}e^{-2\rho\sqrt{\nu{k}}})$. Combining this with
the error bound of $f_n(t)$ given above, the desired error estimate
\eqref{eq:WeeksRate} follows. This ends the proof.
\end{proof}

\begin{figure}[ht]
\centering
\includegraphics[width=0.48\textwidth,height=0.42\textwidth]{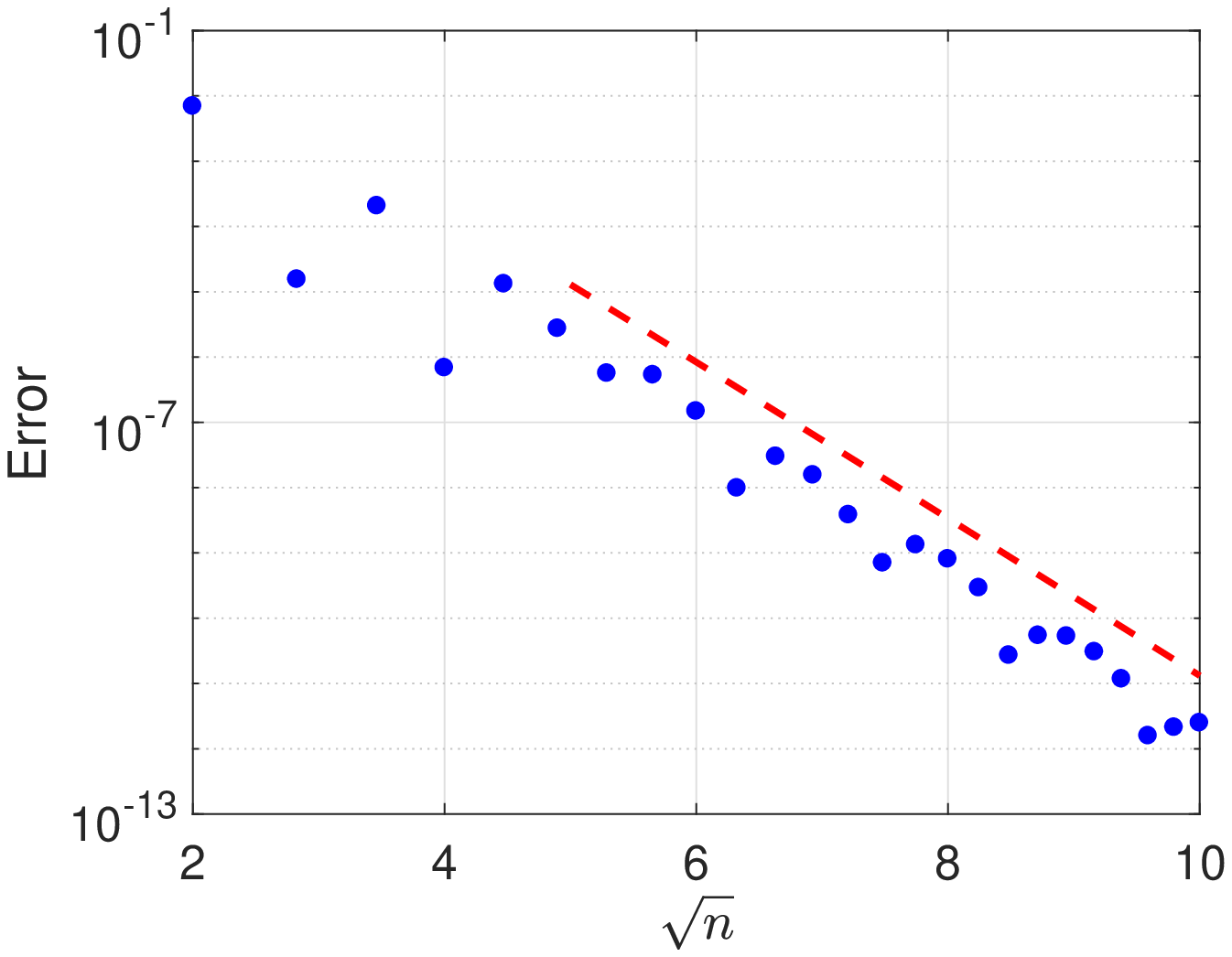}\quad
\includegraphics[width=0.48\textwidth,height=0.42\textwidth]{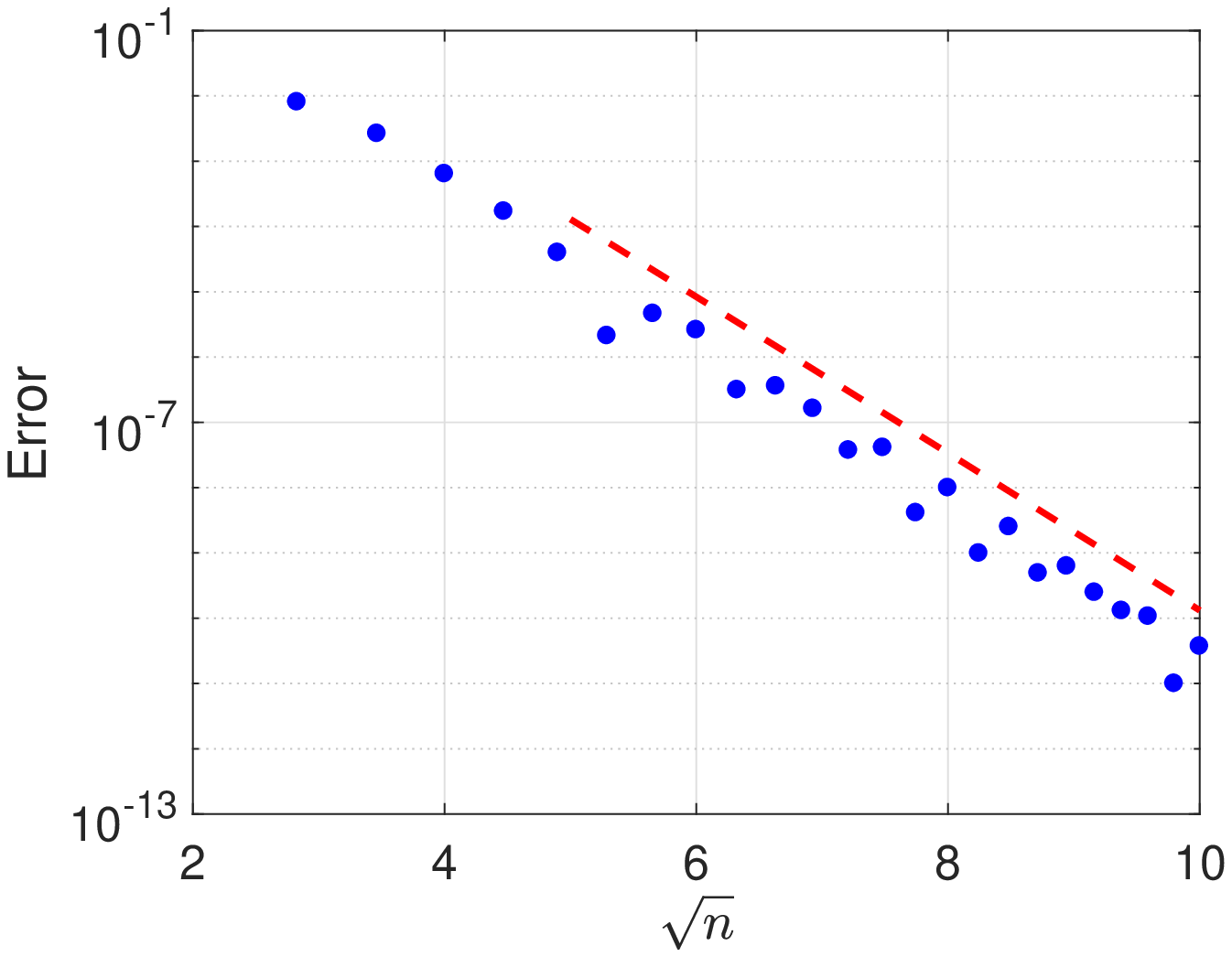}
\caption{The errors of the Weeks method $f_n(t)$ at $t=1$ (left) and
$t=5$ (right) as a function of $\sqrt{n}$. The dash lines show the
predicted rate $O(n^{1/4}\exp(-2\sqrt{2n}))$.}\label{fig:WeekMethod}
\end{figure}

To validate our analysis, we consider the following example
\begin{equation}
f(t) = \frac{1}{t+1}, \quad  (\mathcal{L}f)(s)=e^{s}E_1(s), \notag
\end{equation}
where $E_1(s)$ is the exponential integral, which has been discussed
in \cite{Weideman2023}. When choosing $\sigma=1$ and $\nu=2$, it has
been calculated by empirical fit in \cite{Weideman2023} that the
convergence rate of the Weeks method is approximately
$O(\exp(-C\sqrt{n}))$ with $C\approx2.9$. In fact, since $\rho=1$
for this example, we can conclude from Theorem \ref{thm:Weeks} that
the precise rate of convergence of the Weeks method is
$O(n^{1/4}\exp(-C\sqrt{n}))$ with $C=2\sqrt{2}\approx2.828$. In
Figure \ref{fig:WeekMethod} we illustrate the errors of the Weeks
method at two values of $t$. We see that the actual rates of
convergence of the Weeks method are consistent with our predictions.

\begin{remark}
The error estimate given in Theorem \ref{thm:Weeks} is sharp. To show this, for the example considered above, the Laguerre expansion in \eqref{eq:LagS} becomes
\begin{equation*}
f(t) = \sum_{k=0}^{\infty} c_k L_k(2t), 
\end{equation*}
and thus $\{c_k\}_{k=0}^{\infty}$ are the Laguerre coefficients of $f(t/2)$. From \cite{Elliott1974} we know that $c_k = 2\Gamma(k+1) U(k+1,1;2)$ and it is easily verified that $\{c_k\}_{k=0}^{\infty}$ is a positive and strictly monotonically decreasing sequence. Moreover, by \eqref{def:KummerUAsy},
\[
c_k = 2\sqrt{\pi} \frac{e^{1-2\sqrt{2k}}}{(2k)^{1/4}}(1 + O(k^{-1/2})) = O(k^{-1/4} e^{-2\sqrt{2k}}).
\]
We consider the error of the Weeks method at $t=0$. By the aliasing formula and the monotone decreasing property of $\{c_k\}_{k=0}^{\infty}$, we have
\[
c_k - c_k^w = \tau_k \sum_{\ell=1}^{\infty} (-1)^{\ell+1} \left(
c_{2\ell(n+1)-k} + c_{2\ell(n+1)+k} \right) > 0, \quad k=0,\ldots,n.
\]
Combining these with the fact that $L_k(0)=1$ for $k\in\mathbb{N}_0$, we deduce that
\begin{align}
f(0) - f_n(0) &= \sum_{k=0}^{n} (c_k - c_k^w) + \sum_{k=n+1}^{\infty} c_k > \sum_{k=n+1}^{\infty} c_k = O(n^{1/4} e^{-2\sqrt{2n}}), \nonumber
\end{align}
where we have used the fact that $c_k>0$ and $c_k=O(k^{-1/4}\exp(-2\sqrt{2k}))$ in the last step.
Hence, the predicted convergence rate by Theorem \ref{thm:Weeks} is sharp.
\end{remark}

\subsection{The scaling factor}
\begin{figure}[ht]
\centering
\includegraphics[width=0.48\textwidth,height=0.42\textwidth]{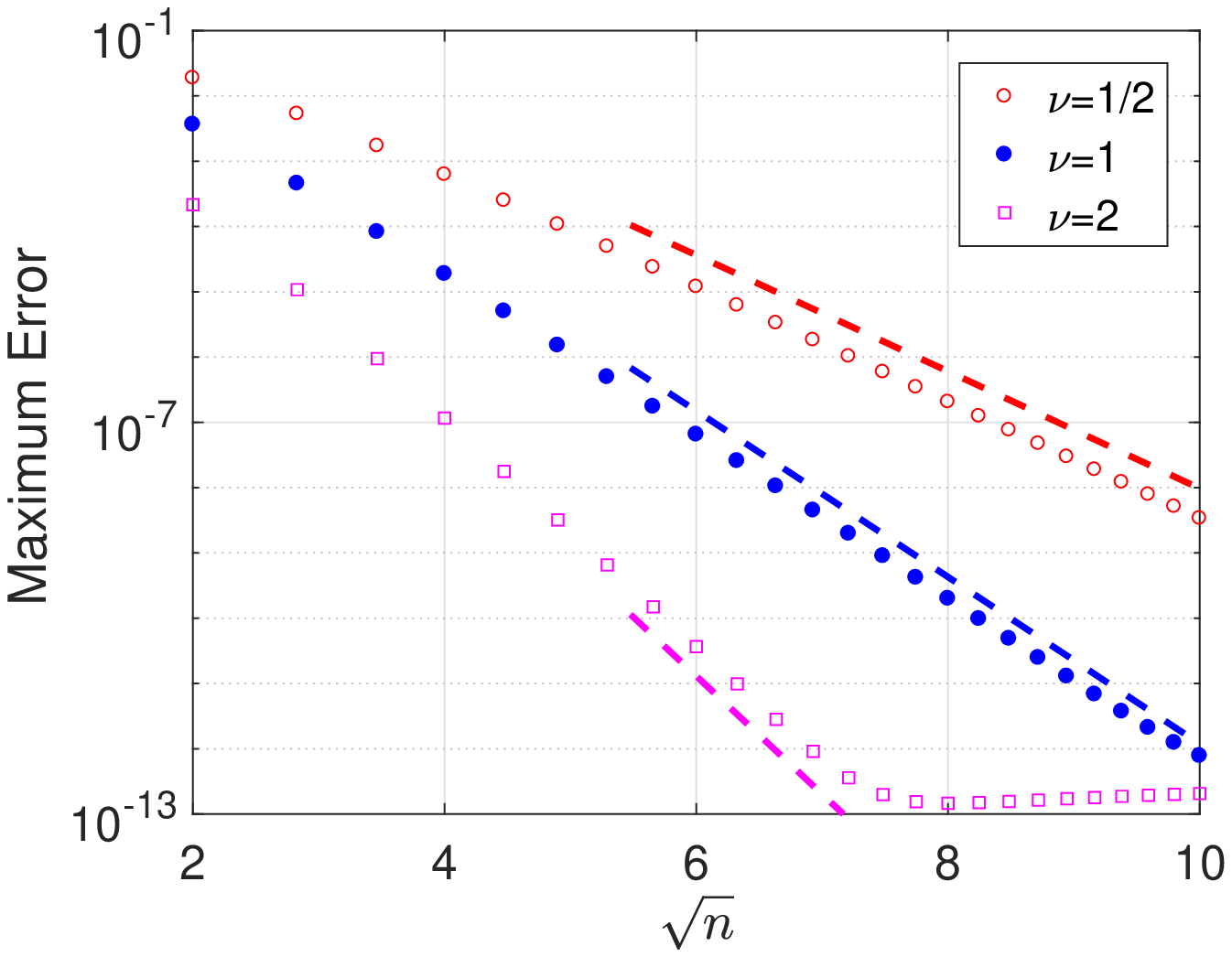}
\includegraphics[width=0.48\textwidth,height=0.42\textwidth]{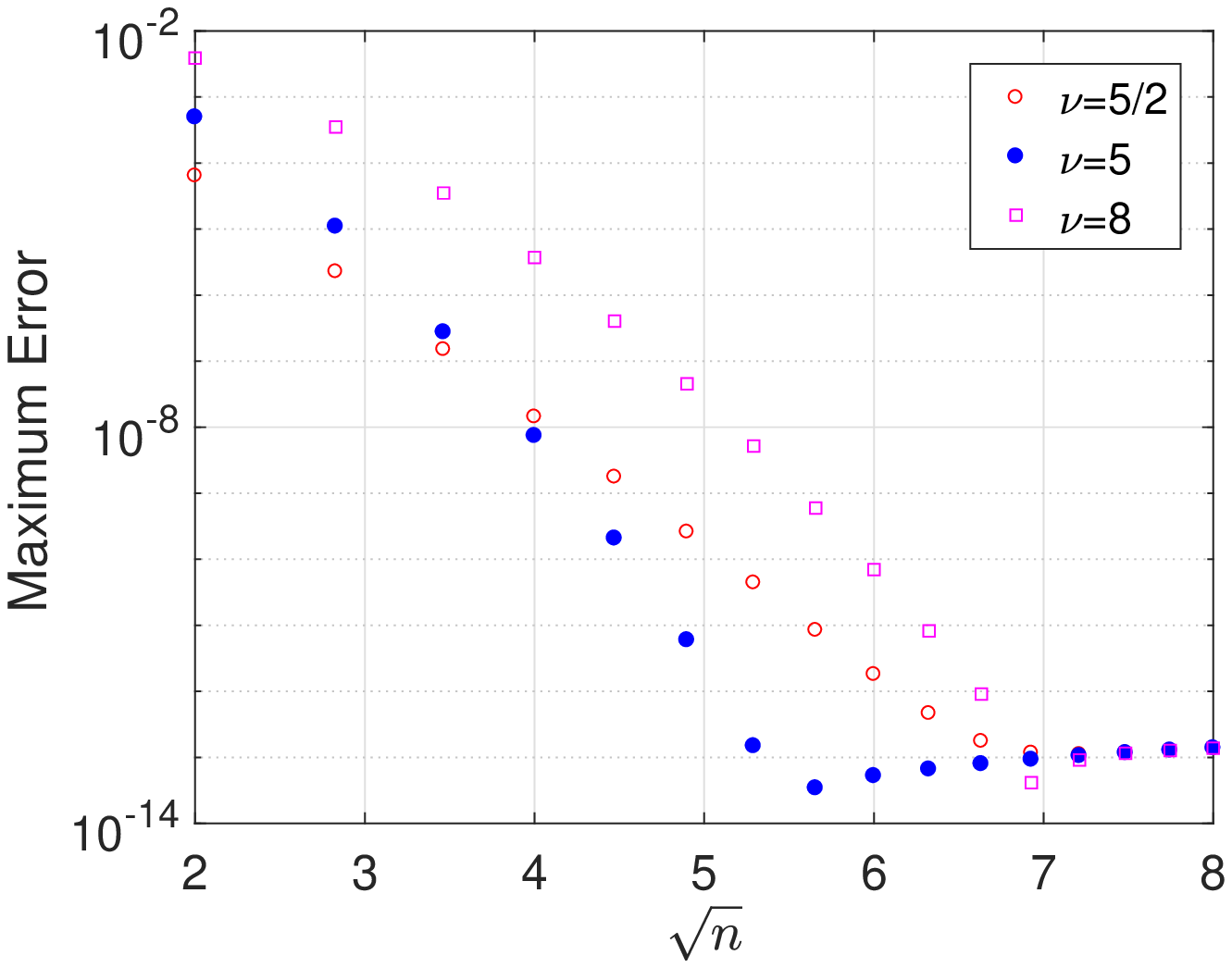}
\caption{The maximum errors of $f_n^{\mathrm{SF}}(x)$ as a function of $\sqrt{n}$ for several different values of $\nu$. The left panel shows $\nu=1/2,1,2$ and the predicted rates $O(n^{1/4}\exp(-3\sqrt{\nu n}))$ (dashed). The right panel shows $\nu=5/2,5,8$. }\label{fig:ScalingFactor}
\end{figure}

In practice, it is often useful to introduce a scaling factor to accelerate the convergence of Laguerre approximations. Specifically, let $f_n^{\mathrm{SF}}(x)$ denote the Laguerre approximation of the form
\begin{equation}
f_n^{\mathrm{SF}}(x) = \sum_{k=0}^{n} c_k \widehat{L}_k^{(\alpha)}(\nu x), \quad c_k = \frac{\nu^{\alpha+1}}{\gamma_k^{(\alpha)}} \int_{0}^{\infty} f(x) \widehat{L}_k^{(\alpha)}(\nu x) \varpi_{\alpha}(x) \mathrm{d}x,
\end{equation}
where $\nu>0$ is a scaling factor. Under the conditions that $f(z)$ is analytic within and on the parabola $P_{\rho}$ for some $\rho>0$ and $|f(z/\nu)|\leq \mathcal{K}|z|^{\beta}e^{-\Re(z)/2}$ for some $\beta\in\mathbb{R}$ as $|z|\rightarrow\infty$ in $D_{\sqrt{\nu}\rho}$, by Theorem \ref{thm:LagFunBound} we can deduce that the convergence rate of $f_n^{\mathrm{SF}}(x)$ in the maximum norm is $O(n^{|\alpha|/2+1/4}e^{-2\rho\sqrt{\nu n}})$. Therefore, once the above conditions on $f$ are satisfied, then the scaling factor $\nu$ should be chosen as large as possible to maximize the convergence rate of $f_n^{\mathrm{SF}}(x)$. To illustrate this, we consider the example $f(x)=e^{-x}/(9+4x)$ and we choose $\alpha=0$ in our computations. By Theorem \ref{thm:LagFunBound}, it is easily verified that the convergence rate of $f_n^{\mathrm{SF}}(x)$ in the maximum norm is $O(n^{1/4}e^{-3\sqrt{\nu n}})$ for $\nu\leq2$. In Figure \ref{fig:ScalingFactor} we plot the maximum errors of $f_n^{\mathrm{SF}}(x)$ for six different values of $\nu$. From the left panel we see that $f_n^{\mathrm{SF}}(x)$ indeed converges at the rate $O(n^{1/4}e^{-3\sqrt{\nu n}})$ for $\nu\leq2$. From the right panel we see that $f_n^{\mathrm{SF}}(x)$ with $\nu=5$ achieves the fastest convergence rate. Moreover, we also see that the convergence rate for $\nu=5$ appears to be faster than root-exponential.

We conclude that our convergence analysis is helpful to some extent in finding scaling factors to accelerate the convergence of Laguerre approximations. However, the problem of finding the optimal scaling factor is still open and we leave it as a topic for future research.

\section{Conclusions}\label{sec:Conclusion}
In this work we have conducted a rigorous convergence rate analysis
of Laguerre approximations of analytic functions. Based on contour
integral techniques, we proved that Laguerre approximations
(including projection and interpolation) converge root-exponentially
when the underlying function is analytic in a parabolic region
containing the positive real axis. We also extended our analysis to
several applications that are of practical interest, including
Laguerre spectral differentiations, Gauss-Laguerre quadrature rules, the scaling factors
and the Weeks method for the inversion of Laplace transform.
In particular, we obtained some sharp estimates for the convergence
rates of Gauss-Laguerre quadrature rules and the Weeks method for
the inversion of Laplace transform of analytic functions.

Before closing this paper, we remark that it is possible to extend
the current study to the convergence rate analysis of Hermite
spectral approximations of analytic functions. Note that this issue
has been studied in the past few decades (see
\cite{Boyd1984,Elliott1974,Hille1940,Shen2009,Xiang2012} and
references therein). However, these studies are partial and
fragmented and further research efforts are required to address this
issue more fully. We will report our results in a forthcoming
publication.

\section*{Acknowledgements}
The author wishes to thank the anonymous referee for his/her valuable comments on this
work. 

\bibliographystyle{amsplain}

\begin{thebibliography}{10}
\bibitem{Barrett1961}
W. Barrett, Convergence properties of Gaussian quadrature formulae,
Comput. J., 3(4):272--277, 1961.

\bibitem{Bern1912}
S. Bernstein, Sur l'ordre de la meilleure approximation des
fonctions continues par les polyn\^{o}mes de degr\'{e} donn\'{e},
Mem. Cl. Sci. Acad. Roy. Belg., 4:1--103, 1912.


\bibitem{Boyd1984}
J. P. Boyd, Asymptotic coefficients of Hermite function series, J.
Comput. Phys., 54(3):382--410, 1984.

\bibitem{Boyd2000}
J. P. Boyd, Chebyshev and Fourier Spectral Methods, Second Edition,
Dover Publications, New York, 2000.

\bibitem{Cohen2007}
A. M. Cohen, Numerical Methods for Laplace transform Inversion,
Numerical Methods and Algorithms, Vol. 5, Springer, New York, 2007.



\bibitem{Davis1975}
P. J. Davis, Interpolation and Approximation, Dover Publications,
New York, 1975.

\bibitem{Davis1984}
P. J. Davis and P. Rabinowitz, Methods of Numerical Integration,
Second Edition, Academic Press, New York, 1984.

\bibitem{Elliott1964}
D. Elliott, The evaluation and estimation of the coefficients in the
Chebyshev series expansion of a function, Math. Comp., 18:274--284,
1964.

\bibitem{Elliott1974}
D. Elliott and P. D. Tuan, Asymptotic estimates of Fourier
coefficients, SIAM J. Math. Anal., 5(1):1--10, 1974.


\bibitem{Hille1940}
E. Hille, Contributions to the theory of Hermitian series II. The
representation problem, Trans. Amer. Math. Soc., 47(1):80--94, 1940.

\bibitem{Giunta1988}
G. Giunta, G. Laccetti and M. R. Rizzardi, More on the Weeks method
for the numerical inversion of the Laplace transform, Numer. Math.,
54:193--200, 1988.

\bibitem{Grad2007}
I. S. Gradshteyn and I. M. Ryzhik, Table of Integrals, Series, and
Products, Seventh Edition, Academic Press, 2007.





\bibitem{Ismail2005}
M. E. H. Ismail, Classical and Quantum Orthogonal Polynomials in One
Variable, Cambridge University Press, Cambridge, 2005.

\bibitem{Liu2019}
W.-J. Liu, L.-L. Wang and H.-Y. Li, Optimal error estimates for
Chebyshev approximation of functions with limited regularity in
fractional Sobolev-type spaces, Math. Comp., 88(320):2857--2895,
2019.




\bibitem{Muskhl1972}
N. I. Muskhelishvili, Singular Integral Equations: Boundary problems
of functions theory and their applications to mathematical physics,
Wolters-Noordhoff Publishing, Groningen, 1972.


\bibitem{Olver2010}
F. W. J. Olver, D. W. Lozier, R. F. Boisvert and C. W. Clark, NIST
Handbook of Mathematical Functions, Cambridge University Press,
Cambridge, 2010.


\bibitem{Opsomer2022}
P. Opsomer and D. Huybrechs, High-order asymptotic expansions of
Gaussian quadrature rules with classical and generalised weight
functions, J. Comput. Appl. Math.,
\texttt{doi:https://doi.org/10.1016/j.cam.2023.115317}, 2023.


\bibitem{Shen2000}
J. Shen, Stable and efficient spectral methods in unbounded domains
using Laguerre functions, SIAM J. Numer. Anal., 38(4):1113--1133,
2000.

\bibitem{Shen2009}
J. Shen and L.-L. Wang, Some recent advances on spectral methods for
unbounded domains, Commun. Comput. Phys., 5(2-4):195--241, 2009.

\bibitem{Shen2011}
J. Shen, T. Tang and L.-L. Wang, Spectral Methods: Algorithms,
Analysis and Applications, Springer, Heidelberg, 2011.


\bibitem{Szasz1958}
O. Sz\'{a}sz and N. Yeardley, The representation of an analytic
function by general Laguerre series, Pacific J. Math., 8:621--633,
1958.

\bibitem{Szego1975}
G. Szeg\H{o}, Orthogonal Polynomials, Fourth Edition, American
Mathematical Society Colloquium Publications, Volume~23, Providence,
Rhode Island, 1975.


\bibitem{Trefeth2013}
L. N. Trefethen, Approximation Theory and Approximation Practice,
SIAM, Philadelphia, 2013.

\bibitem{Wang2012}
H.-Y. Wang and S.-H. Xiang, On the convergence rates of Legendre
approximations, Math. Comp., 81(278):861--877, 2012.

\bibitem{Wang2014}
H.-Y. Wang, D. Huybrechs and S. Vandewalle, Explicit barycentric
weights for polynomial interpolation in the roots or extrema of
classical orthogonal polynomials, Math. Comp., 83(290):2893--2914,
2014.

\bibitem{Wang2016}
H.-Y. Wang, On the optimal estimates and comparison of Gegenbauer
expansion coefficients, SIAM J. Numer. Anal., 54(3):1557--1581,
2016.


\bibitem{Wang2021}
H.-Y. Wang, How much faster does the best polynomial approximation
converge than Legendre projection?, Numer. Math., 147(2):481--503,
2021.

\bibitem{Wang2022}
H.-Y. Wang, Optimal rates of convergence and error localization of
Gegenbauer projections, IMA J. Numer. Anal., 43(4):2413--2444, 2023.

\bibitem{Wang2023}
H.-Y. Wang, Analysis of error localization of Chebyshev spectral
approximations, SIAM J. Numer. Anal., 61(2):952--972, 2023.


\bibitem{Weeks1966}
W. T. Weeks, Numerical inversion of Laplace transforms using
Laguerre functions, J. ACM, 13(3):419--426, 1966.

\bibitem{Weideman1999}
J. A. C. Weideman, Algorithms for parameter selection in the Weeks
method for inverting the Laplace transform, SIAM J. Sci. Comput.,
21(1):111--128, 1999.


\bibitem{Weideman2023}
J. A. C. Weideman and B. Fornberg, Fully numerical Laplace transform
methods, Numer. Algor., 92(1):985--1006, 2023.


\bibitem{Xiang2012}
S.-H. Xiang, Asymptotics on Laguerre or Hermite polynomial
expansions and their applications in Gauss quadrature, J. Math.
Anal. Appl., 393:434--444, 2012.


\bibitem{Xiang2020}
S.-H. Xiang and G.-D. Liu, Optimal decay rates on the asymptotics of
orthogonal polynomial expansions for functions of limited
regularities, Numer. Math., 145:117--148, 2020.

\bibitem{Xiang2021}
S.-H. Xiang, Convergence rates on spectral orthogonal projection
approximation for functions of algebraic and logarithmic
regularities, SIAM J. Numer. Anal., 59(3):1374--1398, 2021.

\bibitem{Xiang2023}
S.-H. Xiang, D.-S. Kong, G.-D. Liu and L.-L. Wang, Pointwise error
estimates and local superconvergence of Jacobi expansions, Math.
Comp., 92:1747-1778, 2023.

\bibitem{Xie2013}
Z.-Q. Xie, L.-L. Wang and X.-D. Zhao, On exponential convergence of
Gegenbauer interpolation and spectral differentiation, Math. Comp.,
82(282):1017--1036, 2013.

\bibitem{Zhao2013}
X.-D. Zhao, L.-L. Wang and Z.-Q. Xie, Sharp error bounds for Jacobi
expansions and Gegenbauer-Gauss quadrature of analytic functions,
SIAM J. Numer. Anal., 51(3):1443--1469, 2013.

\end{thebibliography}

\end{document}